\documentclass{article}

\usepackage{PRIMEarxiv}

\usepackage[utf8]{inputenc} 
\usepackage[T1]{fontenc}    
\usepackage{hyperref}       
\usepackage{url}            
\usepackage{booktabs}       
\usepackage{amsfonts}       
\usepackage{nicefrac}       
\usepackage{microtype}      
\usepackage{lipsum}
\usepackage{fancyhdr}       
\usepackage{graphicx}       
\graphicspath{{media/}}     

\usepackage{amsmath}
\usepackage{amssymb}
\usepackage{epsfig}
\usepackage{epstopdf}
\usepackage{centernot}
\usepackage[algoruled,linesnumbered,noend]{algorithm2e}

\DeclareMathOperator{\Pref}{Pref}

\DeclareMathOperator{\CFL}{CFL}
\DeclareMathOperator{\ICFL}{ICFL}

\newtheorem{proposition}{Proposition}[section]

\newtheorem{definition}{Definition}[section]
\newtheorem{theorem}{Theorem}[section]

\newtheorem{lemma}{Lemma}[section]
\newtheorem{remark}{Remark}[section]
\newtheorem{example}{Example}[section]
\newtheorem{corollary}{Corollary}[section]

\newcommand{\PMCI}{\mathcal{PMC}}

\newcommand{\NB}{\mathcal{NB}}

\numberwithin{equation}{section}

\def\petitcarre{\vrule height4pt width 4pt depth0pt}
\def\enddim{\relax\ifmmode\eqno{\hbox{\petitcarre}}
\else
{\unskip\nobreak\hfil\penalty50
   \hskip2em\hbox{}\nobreak\hfil
   \petitcarre
   \parfillskip=0pt \finalhyphendemerits=0
  \par\medskip}\fi}

\def \begdim {\noindent {\sc Proof} : \par \noindent}

\title{From the Lyndon factorization to the Canonical Inverse Lyndon factorization: back and forth
}

\author{
  Paola Bonizzoni \\
  Dipartimento di Informatica Sistemistica e Comunicazione \\ Universit\`a degli Studi di Milano Bicocca \\
  Milano\\
  \texttt{paola.bonizzoni@unimib.it} \\
   \And
  Clelia De Felice, Rocco Zaccagnino, Rosalba Zizza \\
  Dipartimento di Informatica \\
  Universit\`a degli Studi di Salerno \\
  Salerno\\
  \texttt{\{cdefelice,rzaccagnino,zizza\}@unisa.it} \\
}

\begin{document}
\maketitle

\begin{abstract}
The notion of inverse Lyndon word is related to the classical notion of Lyndon word.
More precisely, inverse Lyndon words are all and only the
nonempty prefixes of the powers of the anti-Lyndon words, where an anti-Lyndon word with respect to a lexicographical order
is a classical Lyndon word with respect to the
inverse lexicographic order.
Each word $w$ admits a factorization in inverse Lyndon words, named the canonical inverse Lyndon factorization
$\ICFL(w)$, which maintains the main properties of the Lyndon factorization of $w$.
Although there is a huge literature on the Lyndon factorization,
the relation between the Lyndon factorization $\CFL_{in}$
with respect to the inverse order and the canonical inverse Lyndon factorization $\ICFL$
has not been thoroughly investigated.
In this paper, we address this question and we
show how to obtain one factorization from the other via the notion of grouping, defined in \cite{inverseLyndon}.
This result naturally opens new insights in the investigation of the relationship between $\ICFL$ and other notions, e.g.,
variants of Burrows Wheeler Transform, as already done for the Lyndon factorization.
\end{abstract}


\section{Introduction}

A word $w$ over a totally ordered alphabet $\Sigma$ is a {\em Lyndon word} if for each
nontrivial factorization $w = uv$, $w$ is strictly smaller than $vu$ for the lexicographical ordering.
A well-known theorem of Lyndon asserts that any nonempty word
factorizes uniquely into a nonincreasing product of Lyndon words, called its
Lyndon factorization \cite{Lyndon}.
It can be efficiently computed.
Linear-time algorithms for computing this factorization
can be found in \cite{duval} whereas an $\mathcal{O}(\lg{n})$-time parallel algorithm
has been proposed in \cite{apostolico-crochemore-1995}.
There are several results which give relations between
Lyndon words, codes and combinatorics of words
and there are many algorithmic applications of the Lyndon factorization \cite{PerRes}.
For instance, bijective versions of the classical Burrows Wheeler Transform
have been defined in \cite{BKKP19,GilScot,varBWT} and they are based on combinatorial results on Lyndon words
proved in \cite{GeReu} (see \cite{GeResReu} for more recent related results).
As another example, a main property proved in \cite{restivo-sorting-2014} is the {\em Compatibility Property}.
Roughly, the compatibility property allows us to extend the mutual order between suffixes
of products of the Lyndon factors to the suffixes of the whole word.
Similar ideas have been further investigated to accelerate suffix sorting in practice \cite{lyndonAcc}.

More recently Lyndon words found a renewed
theoretical interest and several variants of them have been studied \cite{Nyldon2,DRR,PZ20,BW20,NyldonNew}.
In particular, inverse Lyndon words have been introduced in \cite{inverseLyndon}.
More precisely, an anti-Lyndon word with respect to a lexicographical order
is a classical Lyndon word with respect to the
inverse lexicographic order. Then, inverse Lyndon words are all and only the
nonempty prefixes of the powers of the anti-Lyndon words.
Each word $w$ admits a factorization in inverse Lyndon words, named the canonical inverse Lyndon factorization
$\ICFL(w)$, which maintains the main properties of the Lyndon factorization of $w$: it is uniquely determined and can be computed in linear time.
In addition, it maintains a similar Compatibility Property.
Other combinatorial properties of $\ICFL(w)$ have been proved in \cite{inverseLyndonTCS21}.

In this paper we address a main theoretical question left open in \cite{inverseLyndon}.
More precisely, we have proved in \cite{inverseLyndon} that the canonical inverse Lyndon factorization
$\ICFL(w)$ of $w$ is a grouping of the Lyndon factorization $\CFL_{in}(w)$
of $w$ with respect to the inverse order.
Roughly, this means that every element in $\ICFL(w)$ is the concatenation
of factors of $\CFL_{in}(w)$ that are related by the prefix order.
The proof of this result is not constructive and leaves open the problem of constructing this grouping,
that is, how to obtain $\ICFL$ from $\CFL_{in}$.
Moreover, another open problem is how to obtain $\CFL_{in}$ from $\ICFL$.
The main result of this paper is to solve these two problems, by providing a full characterization of the
relationship between $\CFL_{in}$ and $\ICFL$ via grouping.

In detail, we know that every element in $\ICFL(w)$ is the concatenation
of factors of $\CFL_{in}(w)$ that are related by the prefix order.
Then, in Proposition \ref{step1} we prove that $\ICFL(w)$
can be computed locally in the non-increasing maximal chains for the prefix order
of $\CFL_{in}(w)$.
There is a pair of words, the canonical pair associated with $w$,
which is crucial in the construction of $\ICFL(w)$ (see Section \ref{icfl} for the
definition).
In Proposition \ref{Algo3} we show how to determine such a pair
in one of the aforementioned chains.

The procedure for obtaining $\CFL_{in}(w)$ from $\ICFL(w)$ is very simple.
We know that every element in $\ICFL(w)$ is the concatenation
of factors of $\CFL_{in}(w)$ that are related by the prefix order.
Therefore, we first demonstrate that we can limit ourselves to handling the single factor
$m_i$ in $\ICFL(w) = (m_1, \ldots , m_k)$ (Corollary \ref{cor2}).
Next, we show that each bordered word has a unique
nonempty border which is unbordered (Proposition \ref{Prdopo4}).
This combinatorial property allows us to recursively construct a sequence of words
uniquely associated with each $m_i$ which will turn out to be $\CFL_{in}(m_i)$
(Definition \ref{def:NB}, Corollary \ref{cor3}).

Propositions \ref{InvPreN} and \ref{bordo1} were proved in
\cite{inverseLyndon} and \cite{inverseLyndonTCS21}, respectively.
Both make use of Proposition \ref{preprime} which was stated in a slightly incorrect form
in \cite{inverseLyndon,inverseLyndonTCS21}.
After having correctly stated Proposition \ref{preprime}, as minor results we demonstrate
Propositions \ref{InvPreN} and \ref{bordo1} again.

The paper is organized as follows.
In Sections \ref{W}, \ref{LW}, \ref{LYF}, \ref{ILW}, \ref{icfl}, \ref{Group}, we gathered
the basic definitions and known results we need.
We show how to obtain the Lyndon factorization $\CFL_{in}(w)$
of $w$ with respect to the inverse order from the canonical inverse Lyndon factorization $\ICFL(w)$
in Section \ref{fromItoC}
and we illustrate how to group factors of $\CFL_{in}(w)$ to obtain $\ICFL(w)$
in Section \ref{fromCtoI}.
We conclude in Section \ref{Conc} with some main issues that follow on from
the results proved in the paper.

\section{Words} \label{W}

Throughout this paper we follow
\cite{bpr,CK,Lo,lothaire,reu} for the notations.
We denote by $\Sigma^{*}$ the {\it free monoid}
generated by a finite alphabet $\Sigma$
and we set $\Sigma^+=\Sigma^{*} \setminus 1$, where $1$ is
the empty word.
For a word $w \in \Sigma^*$, we denote by $|w|$ its {\it length}.
A word $x \in \Sigma^*$ is a {\it factor} of $w \in \Sigma^*$ if there are
$u_1,u_2 \in \Sigma^*$ such that $w=u_1xu_2$.
If $u_1 = 1$ (resp. $u_2 = 1$), then $x$ is a {\it prefix}
(resp. {\it suffix}) of $w$.
A factor (resp. prefix, suffix) $x$ of $w$
is {\it proper} if $x \not = w$.
Two words $x,y$ are {\it incomparable} for the prefix order, denoted as $x \Join y$,
if neither $x$ is a prefix of $y$ nor $y$ is a prefix of $x$.
Otherwise, $x,y$ are {\it comparable} for the prefix order.
We write $x \leq_{p} y$ if $x$ is a prefix of $y$
and  $x \geq_{p} y$ if $y$ is a prefix of $x$.
The notion of a pair of words comparable (or incomparable) for the suffix order
is defined symmetrically.

We recall that, given a nonempty word $w$,
a {\it border} of $w$ is a word which is both a proper prefix and a suffix
of $w$ \cite{CHL07}. The longest proper prefix of $w$ which is a suffix of
$w$ is also called {\it the border} of $w$ \cite{CHL07,lothaire}.
A word $w \in \Sigma^+$ is {\em bordered} if it has a nonempty border.
Otherwise, $w$ is {\it unbordered}.
A nonempty word $w$ is \textit{primitive} if
$w = x^k$ implies $k = 1$. An unbordered word is primitive.
A {\it sesquipower} of a word $x$ is a word $w = x^np$ where
$p$ is a proper prefix of $x$ and $n \geq 1$.

Two words $x,y$ are called {\em conjugate} if there exist words
$u,v$ such that $x=uv, y=vu$.
The conjugacy relation is an equivalence relation. A conjugacy class
is a class of this equivalence relation.

\begin{definition} \label{lex-order}
Let $(\Sigma, <)$ be a totally ordered alphabet.
The {\it lexicographic} (or {\it alphabetic order})
$\prec$ on $(\Sigma^*, <)$ is defined by setting $x \prec y$ if
\begin{itemize}
\item $x$ is a proper prefix of $y$,
or
\item $x = ras$, $y =rbt$, $a < b$, for $a,b \in \Sigma$ and $r,s,t \in \Sigma^*$.
\end{itemize}
\end{definition}

In the next part of the paper we will implicitly refer to
totally ordered alphabets.
For two nonempty words $x,y$, we write $x \ll y$ if
$x \prec y$ and $x$ is not a proper prefix of $y$
\cite{Bannai15}. We also write $y \succ x$ if $x \prec y$.
Basic properties of the lexicographic order are recalled below.

\begin{lemma} \label{proplexord}
For $x,y \in \Sigma^*$, the following properties hold.
\begin{itemize}
\item[(1)]
$x \prec y$ if and only if $zx \prec zy$,
for every word $z$.
\item[(2)]
If $x \ll y$, then $xu \ll yv$
for all words $u,v$.
\item[(3)]
If $x \prec y \prec xz$ for a word $z$,
then $y = xy'$ for some word $y'$ such that $y' \prec z$.
\end{itemize}
\end{lemma}

Let $\mathcal{S}_1, \ldots , \mathcal{S}_t$ be sequences, with
$\mathcal{S}_j = (s_{j,1}, \ldots , s_{j,r_j})$.
For abbreviation, we let
$$(\mathcal{S}_1, \ldots , \mathcal{S}_t)$$
stand for the sequence
$$(s_{1,1}, \ldots , s_{1,r_1}, \ldots , s_{t,1}, \ldots , s_{t,r_t})$$


\section{Lyndon words} \label{LW}

\begin{definition}\label{Lyndon-word}
A Lyndon word $w \in \Sigma^+$ is a word which is primitive and the smallest one
in its conjugacy class for the lexicographic order.
\end{definition}

\begin{example}
{\rm Let $\Sigma = \{a,b\}$ with $a < b$.
The words $a$, $b$, $aaab$, $abbb$, $aabab$ and $aababaabb$
are Lyndon words. On the contrary, $abab$, $aba$ and
$abaab$ are not Lyndon words. Indeed, $abab$ is a non-primitive word,
$aab \prec aba$
and $aabab \prec abaab$.}
\end{example}

Lyndon words are also called {\it prime words}
and their prefixes are also called {\it preprime words} in \cite{Knuth}.

\begin{proposition} \label{Pr1}
Each Lyndon word $w$ is unbordered.
\end{proposition}

A class of conjugacy is also called a {\it necklace} and often identified
with the minimal word for the lexicographic
order in it. We will adopt this terminology. Then
a word is a necklace if and only if it is a power of a Lyndon word.
A {\it prenecklace} is a prefix of a necklace. Then clearly
any nonempty prenecklace $w$ has the form $w = (uv)^ku$, where $uv$ is a Lyndon word,
$u \in \Sigma^*$, $v \in \Sigma^+$, $k \geq 1$, that is, $w$ is a sesquipower of a Lyndon word $uv$.
The following result has been proved in \cite{duval}.
It shows that the nonempty prefixes of Lyndon words are exactly
the nonempty prefixes of the powers of Lyndon words with the exclusion of
$c^k$, where $c$ is the maximal letter and $k \geq 2$.

\begin{proposition} \label{preprime}
A word is a nonempty preprime word if and only if it is a sesquipower of a Lyndon word distinct of
$c^k$, where $c$ is the maximal letter and $k \geq 2$.
\end{proposition}

The proof of Proposition \ref{preprime} uses the following result which
characterizes, for a given nonempty prenecklace $w$ and a letter $b$, whether $wb$ is still
a prenecklace or not and, in the first case, whether $wb$ is a Lyndon word or not \cite[Lemma 1.6] {duval}.

\begin{theorem} \label{FundamentalPreneck}
Let $w = (uav')^ku$ be a nonempty prenecklace, where $uav'$ is a Lyndon word,
$u, v' \in \Sigma^*$, $a \in \Sigma$, $k \geq 1$. For any $b \in \Sigma$, the
word $wb$ is a prenecklace if and only if $b \geq a$. Moreover $wb \in L$ if and only if
$b > a$.
\end{theorem}


\section{Lyndon factorization} \label{LYF}

In the following $L = L_{(\Sigma^*, <)}$ will be the set of Lyndon words, totally
ordered by the relation $\prec$ on $(\Sigma^*, <)$.

\begin{theorem} \label{Lyndon-factorization}
Any word $w \in \Sigma^+$ can be written in a unique way as
a nonincreasing product $w=\ell_1 \ell_2 \cdots \ell_h$ of Lyndon words, i.e., in the form
\begin{eqnarray} \label{LF}
w & = & \ell_1 \ell_2 \cdots \ell_h, \mbox{ with } \ell_j \in L \mbox{ and } \ell_1 \succeq \ell_2 \succeq \ldots \succeq \ell_h
\end{eqnarray}
\end{theorem}

The sequence $\CFL(w) = (\ell_1, \ldots, \ell_h)$ in Eq. \eqref{LF} is called the
\textit{Lyndon decomposition} (or \textit{Lyndon factorization}) of $w$.
It is denoted by $\CFL(w)$ because Theorem \ref{Lyndon-factorization}
is usually credited to Chen, Fox and Lyndon
\cite{Lyndon}.
Uniqueness of the above factorization can be obtained as a consequence of the following result, proved
in \cite{duval}.

\begin{lemma} \label{duval-prop}
Let $w \in \Sigma^+$ and let $\CFL(w) = (\ell_1, \ldots, \ell_h)$.
Then the following properties hold:
\begin{itemize}
\item[(i)]
$\ell_h$ is the nonempty suffix of $w$ which is the smallest with respect to the lexicographic order.
\item[(ii)]
$\ell_h$ is the longest
suffix of $w$ which is a Lyndon word.
\item[(iii)]
$\ell_1$ is the longest
prefix of $w$ which is a Lyndon word.
\end{itemize}
\end{lemma}

A direct consequence is stated below and it is necessary for our aims.

\begin{corollary} \label{cor1}
Let $w \in \Sigma^+$, let $\ell_1$ be its longest
prefix which is a Lyndon word and let $w'$ be
such that $w= \ell_1 w'$. If $w' \not = 1$, then $\CFL(w) = (\ell_1, \CFL(w'))$.
\end{corollary}

Sometimes we need to emphasize consecutive equal factors
in $\CFL$. We write
$\CFL(w) = (\ell_1^{n_1}, \ldots, \ell_r^{n_r})$ to denote
a tuple of $n_1 + \ldots + n_r$ Lyndon words,
where $r  > 0$, $n_1, \ldots , n_r \geq 1$. Precisely
$\ell_1 \succ \ldots \succ \ell_r$
are Lyndon words, also named {\it Lyndon factors} of $w$.
There is a linear time algorithm
to compute the pair $(\ell_1, n_1)$ and thus, by iteration,
the Lyndon factorization of $w$ \cite{lothaire,FM}.
Linear time algorithms may also be found in \cite{duval}
and in the more recent paper
\cite{GGT}.


\section{Inverse Lyndon words} \label{ILW}

For the material in this section see \cite{inverseLyndon,lata2020,inverseLyndonTCS21}.
Inverse Lyndon words are related to the
inverse alphabetic order. Its definition is recalled below.

\begin{definition} \label{ILO}
Let $(\Sigma, <)$ be a totally ordered alphabet.
The {\rm inverse} $<_{in}$ of $<$ is defined
by
$$ \forall a, b \in \Sigma \quad b <_{in} a \Leftrightarrow a < b $$
The {\rm inverse lexicographic} or {\rm inverse alphabetic order}
on $(\Sigma^*, <)$, denoted $\prec_{in}$, is the lexicographic order
on $(\Sigma^*, <_{in})$.
\end{definition}

\begin{example}
{\rm Let $\Sigma = \{a,b,c,d \}$ with $a < b < c < d$. Then $dab \prec dabd$
and $dabda \prec dac$. We have
$d <_{in} c <_{in} b <_{in} a$. Therefore $dab \prec_{in} dabd$ and $dac \prec_{in} dabda$.}
\end{example}

Of course for all $x, y \in \Sigma^*$ such that $x \Join y$,
$$y \prec_{in} x \Leftrightarrow x \prec y .$$
Moreover, in this case $x \ll y$. This justifies
the adopted terminology.

From now on, $L_{in} = L_{(\Sigma^*, <_{in})}$ denotes the set
of the Lyndon words on $\Sigma^*$ with respect to the inverse lexicographic order.
Following \cite{antiLyndon},
a word $w \in L_{in}$ will be named an {\it anti-Lyndon word}. Correspondingly, an
{\it anti-prenecklace} will be a prefix of an {\it anti-necklace}, which in turn will
be a necklace with respect to the inverse lexicographic order.

In the following, we denote by $\CFL_{in}(w)$ the Lyndon factorization
of $w$ with respect to the inverse order $<_{in}$.

\begin{definition} \label{inverse-Lyndon-word}
A word $w \in \Sigma^+$ is an inverse Lyndon word if
$s \prec w$, for each nonempty proper suffix $s$ of $w$.
\end{definition}

\begin{example}
{\rm The words $a$, $b$, $aaaaa$, $bbba$, $baaab$, $bbaba$ and $bbababbaa$
are inverse Lyndon words on $\{a,b\}$, with $a < b$.
On the contrary, $aaba$ is not an inverse Lyndon word since $aaba \prec ba$.
Analogously, $aabba \prec ba$ and thus $aabba$ is not an inverse Lyndon word.}
\end{example}

The following result has been stated in \cite{inverseLyndonTCS21}.

\begin{proposition} \label{Pr2}
A word $w \in \Sigma^+$ is an anti-Lyndon word if and only if it is
an unbordered inverse Lyndon word.
\end{proposition}

The following results have been proved in \cite{inverseLyndon}.

\begin{lemma} \label{lem:inverse-Lyndon-word-prefix}
Any nonempty prefix of an
inverse Lyndon word is an inverse Lyndon word.
\end{lemma}

The proof of Proposition \ref{InvPreN} was given in \cite{inverseLyndon}.
The reason we report it here is that this proof uses the statement in
Proposition \ref{preprime} but in \cite{inverseLyndon} that statement was slightly incorrect.
However, the proof still holds.

\begin{proposition} \label{InvPreN}
A word $w \in \Sigma^+$ is an inverse Lyndon word if and only if
$w$ is a nonempty anti-prenecklace.
\end{proposition}
\begdim
Let $w \in \Sigma^+$ be an inverse Lyndon word. The first letter of $w$ is an anti-Lyndon word thus
also a nonempty anti-prenecklace.
Let $p$ be the longest nonempty prefix of $w$ which is an anti-prenecklace.
By Theorem \ref{FundamentalPreneck}, if $p$ were distinct from $w$, then we would have
$p = (uav')^ku$, $w = (uav')^kubt$, where $uav' \in L_{in}$,
$u, v', t \in \Sigma^*$, $a, b \in \Sigma$, $a < b$, $k \geq 1$. Thus,
$w \ll ubt$, in contradiction with Definition \ref{inverse-Lyndon-word}.
Therefore, $w$ is a nonempty anti-prenecklace.
Conversely, let $w$ be a nonempty anti-prenecklace, that is, a sesquipower of an anti-Lyndon word.
If $w = a^n$, where $a$ is the minimal letter in $\Sigma$ and $n \geq 2$,
then clearly $w$ is an inverse Lyndon word. Otherwise, by Proposition \ref{preprime},
there is a word $t$ such that $wt \in L_{in}$. By Proposition \ref{Pr2},
$wt$ is an inverse Lyndon word.
If there existed a
nonempty proper suffix $s$ of $w$ such that $w \prec s$, we clearly would have $w \ll s$.
Hence, by item (2) in Lemma \ref{proplexord}, $wt \ll st$, where $st$ is a nonempty proper suffix
of $wt$. This is in contradiction with Definition \ref{inverse-Lyndon-word}, thus $w$ is an inverse Lyndon word.
\enddim


\section{Inverse Lyndon factorizations} \label{icfl}

For the material in this section see \cite{inverseLyndon,lata2020,inverseLyndonTCS21}.
An inverse Lyndon factorization of
a word $w \in \Sigma^+$ is a sequence $(m_1, \ldots, m_{k})$ of inverse Lyndon words
such that $m_1 \cdots m_{k} = w$ and $m_i \ll m_{i+1}$, $1 \leq i \leq k-1$.
A word may have different
inverse Lyndon factorizations (see Example \ref{nonunique})
but it has a unique canonical inverse Lyndon factorization, denoted $\ICFL(w)$.
If $w$ is an inverse Lyndon word, then $\ICFL(w) = w$.
Otherwise, $\ICFL(w)$ is recursively defined.
The first factor of $\ICFL(w)$
is obtained by a special factorization of the shortest nonempty
prefix $z$ of $w$ such that $z$ is not an inverse Lyndon word
defined below.

\begin{definition} \label{brepref} \cite{inverseLyndon}
Let $w \in \Sigma^+$, let $p$ be an inverse Lyndon word
which is a nonempty proper prefix
of $w = pv$.
The {\rm bounded right extension} $\overline{p}$ of $p$
(relatively to $w$), if it exists,
is a nonempty prefix of $v$ such that:
\begin{itemize}
\item[(1)]
$\overline{p}$ is an inverse Lyndon word,
\item[(2)]
$pz'$ is an inverse Lyndon word, for each proper nonempty prefix $z'$
of $\overline{p}$,
\item[(3)]
$p \overline{p}$ is not an inverse Lyndon word,
\item[(4)]
$p \ll \overline{p}$.
\end{itemize}
Moreover, we set
$\Pref_{bre}(w)= \{(p,\overline{p}) ~|~ p \mbox{ is an inverse Lyndon word which is a} \\
\mbox{nonempty proper prefix of } w \}.$
\end{definition}

It has been proved that $\Pref_{bre}(w)$ is empty if and only if
$w$ is an inverse Lyndon word (Proposition 4.2 in \cite{inverseLyndon}).
If $w$ is not an inverse Lyndon word, then
$\Pref_{bre}(w)$ contains only one pair
and the description of this pair
is given below (Propositions 4.1 and 4.3 in \cite{inverseLyndon}).
In the next, the unique pair $(p,\overline{p})$ in $\Pref_{bre}(w)$
will be named {\it the canonical pair} associated with $w$.

\begin{proposition} \label{shortest}
Let $w \in \Sigma^+$ be a word which is not an inverse Lyndon word.
Let $z$ be the shortest nonempty prefix of $w$ which is not an inverse Lyndon word. Then,
\begin{itemize}
\item[(1)]
$z = p \overline{p}$, with $(p , \overline{p}) \in \Pref_{bre}(w)$.
\item[(2)]
$p = ras$ and $\overline{p} = rb$, where $r,s \in \Sigma^*$, $a,b \in \Sigma$ and
$r$ is the shortest prefix of $p \overline{p}$ such that $p \overline{p} = rasrb$,
with $a < b$.
\end{itemize}
\end{proposition}

\begin{example} \label{LATAbre}
{\rm Let $\Sigma = \{a,b \}$ with $a < b$.
Let us consider $w = babaaabb$
and the prefixes $p_1 = bab$ and $p_2 = babaaa$ of $w$.
First, $w$ is not an inverse Lyndon word.
Thus, $\Pref_{bre}(w)$ contains only one pair.
Moreover
each proper nonempty prefix of $w$ is an inverse Lyndon word.
By item (1) in Proposition \ref{shortest}, we have
$w = p \overline{p}$. By item (2) in Proposition \ref{shortest},
the bounded right extension
of $p_1 = bab$ does not exist (we should have $\overline{p_1} = aaabb$ in
contradiction with $p_1 \ll \overline{p_1}$).
Since $w$ starts with $b$, the shortest common prefix $r$ of $p$ and $\overline{p}$
has a positive length. Indeed, $p = p_2 = babaaa$ and
$\overline{p} = \overline{p_2} = bb$.}
\end{example}

The above results suggest the following characterization of the canonical
pair $(p , \overline{p})$ associated with $w$ (Proposition 6.2 in \cite{inverseLyndonTCS21}).

\begin{proposition} \label{characterization}
Let $w \in \Sigma^+$ be a word which is not an inverse Lyndon word.
A pair of words $(p , \overline{p})$ is the canonical pair associated with $w$
if and only the following conditions are satisfied.
\begin{itemize}
\item[(1)]
$z = p \overline{p}$ is the shortest nonempty prefix of $w$ which is not an inverse Lyndon word.
\item[(2)]
$p = ras$ and $\overline{p} = rb$, where $r,s \in \Sigma^*$, $a,b \in \Sigma$ and
$r$ is the shortest prefix of $p \overline{p}$ such that $p \overline{p} = rasrb$,
with $a < b$.
\item[(3)]
$\overline{p}$ is an inverse Lyndon word.
\end{itemize}
\end{proposition}

Given a word $w$ which is not an inverse Lyndon word,
Proposition \ref{characterization} suggests a method to identify
the canonical pair $(p , \overline{p})$ associated with $w$:
just find the shortest nonempty prefix $z$ of $w$ which is not an inverse Lyndon word
and then a factorization $z = p \overline{p}$ such that conditions (2) and (3) in Proposition
\ref{characterization} are satisfied.

The canonical inverse Lyndon factorization has been also recursively defined.

\begin{definition} \label{def:ICFL}
Let $w \in \Sigma^+$. \\
{\rm (Basis Step)}
If $w$ is an inverse Lyndon word,
then $\ICFL(w) = (w)$. \\
{\rm (Recursive Step)}
If $w$ is not an inverse Lyndon word,
let $(p,\overline{p})$ be the canonical pair associated with $w$ and let
$v \in \Sigma^*$ such that $w = pv$.
Let $\ICFL(v) = (m'_1, \ldots, m'_{k})$ and let
$r,s \in \Sigma^*$, $a,b \in \Sigma$ such that $p = ras$, $\overline{p} = rb$ with $a < b$.
$$\ICFL(w) = \begin{cases} (p, \ICFL(v)) & \mbox{ if } \overline{p} = rb \leq_{p} m'_1 \\
(p m'_1, m'_2, \ldots, m'_{k}) & \mbox{ if } m'_1 \leq_{p} r \end{cases}$$
\end{definition}

\begin{example} \label{nonunique} \cite{inverseLyndon}.
{\rm Let $\Sigma = \{a,b,c,d \}$ with $a < b < c < d$,
$w = dabadabdabdadac$. We have
$\CFL_{in}(w) = (daba, dab, dab, dadac)$,
$\ICFL(w) = (daba, dabdab, dadac)$.
Another inverse Lyndon factorizations of $w$
is $(dabadab, dabda, dac)$.
Consider
$z = dabdadacddbdc$.
It is easy to see that
$(dab,dadacd,db,dc)$, $(dabda,dac,ddbdc)$, $(dab, dadac, ddbdc)$
are all inverse Lyndon factorizations of
$z$. The first factorization has four factors whereas the
others have three factors.
Moreover $\ICFL(z) = \CFL_{in}(z) = (dab,dadac,ddbdc)$.}
\end{example}


\section{Groupings} \label{Group}

In this section we recall a special property of $\ICFL$ proved in
\cite{inverseLyndon}, needed here for our aims.

Let $\CFL_{in}(w) = (\ell_1, \ldots , \ell_h)$, where
$\ell_1 \succeq_{in} \ell_2 \succeq_{in} \ldots \succeq_{in} \ell_h$.
Consider the partial order $\geq_p$, where
$x \geq_p y$ if $y$ is a prefix of $x$. Recall that a {\it chain} is a set
of a pairwise comparable elements. We say that a chain is maximal if it is
not strictly contained in any other chain.
A non-increasing {\it (maximal) chain}
in $\CFL_{in}(w)$ is the
sequence corresponding to a (maximal) chain in the
multiset $\{\ell_1, \ldots , \ell_h \}$ with respect to $\geq_p$.
We denote by $\PMCI$ a non-increasing maximal chain in $\CFL_{in}(w)$.
Looking at the definition of the (inverse) lexicographic order,
it is easy to see that a $\PMCI$ is a sequence
of consecutive factors in $\CFL_{in}(w)$.
Moreover $\CFL_{in}(w)$ is
the concatenation of its $\PMCI$.
The formal definitions are given below.

\begin{definition} \label{MaxCh}
Let $w \in \Sigma^+$, let $\CFL_{in}(w) = (\ell_1, \ldots , \ell_h)$
and let $1 \leq r < s \leq h$.
We say that $\ell_{r}, \ell_{r+1}, \ldots , \ell_{s}$
is a non-increasing {\rm maximal chain for the prefix order}
in $\CFL_{in}(w)$, abbreviated $\PMCI$, if
$\ell_{r} \geq_p \ell_{r+1} \geq_p \ldots  \geq_p \ell_{s}$.
Moreover, if $r > 1$, then $\ell_{r - 1} \not \geq_p \ell_{r}$,
if $s < h$, then $\ell_{s} \not \geq_p \ell_{s +1}$.
Two $\PMCI$ $\mathcal{C}_1 = \ell_{r}, \ell_{r+1}, \ldots , \ell_{s}$,
$\mathcal{C}_2 = \ell_{r'}, \ell_{r'+1}, \ldots , \ell_{s'}$ are
{\rm consecutive} if $r' = s+1$ (or $r = s' +1$).
\end{definition}

\begin{definition} \label{DecMaxCh}
Let $w \in \Sigma^+$, let $\CFL_{in}(w) = (\ell_1, \ldots , \ell_h)$.
We say that $(\mathcal{C}_{1}, \mathcal{C}_{2}, \ldots , \mathcal{C}_{s})$
is {\rm the decomposition} of $\CFL_{in}(w)$ into its non-increasing maximal chains for the prefix order
if the following holds
\begin{itemize}
\item[(1)]
Each $\mathcal{C}_{j}$ is a non-increasing maximal chain in $\CFL_{in}(w)$.
\item[(2)]
$\mathcal{C}_{j}$ and $\mathcal{C}_{j + 1}$ are consecutive, $1 \leq j \leq s-1$.
\item[(3)]
$\CFL_{in}(w)$ is the concatenation of the sequences $\mathcal{C}_{1}, \mathcal{C}_{2}, \ldots , \mathcal{C}_{s}$.
\end{itemize}
\end{definition}

\begin{example} \label{maximalchain} \cite{inverseLyndon}.
{\rm Let $\Sigma = \{a,b,c,d \}$ with $a < b < c < d$,
$w = dabadabdabdadac$.
In Example \ref{nonunique}, we observed that
$\CFL_{in}(w) = (daba, dab, dab, dadac)$.
This sequence has two $\PMCI$, namely
$(daba, dab, dab)$, $(dadac)$.
The decomposition of $\CFL_{in}(w)$ into its $\PMCI$
is $((daba, dab, dab), (dadac))$.
Let $z = dabdadacddbdc$. Then
$\CFL_{in}(z) = (dab,dadac,ddbdc)$ has
three $\PMCI$: $(dab)$, $(dadac)$, $(ddbdc)$.
The decomposition of $\CFL_{in}(w)$ into its $\PMCI$
is $((dab), (dadac), (ddbdc))$.}
\end{example}

A {\it grouping} of $\CFL_{in}(w)$ is an inverse Lyndon factorization
$(m_1, \ldots, m_{k})$ of $w$
such that any factor is a product of consecutive factors in a
$\PMCI$ of $\CFL_{in}(w)$. Formally, the
definition of a grouping of $\CFL_{in}(w)$ is given below in two steps. We first
define the grouping of a $\PMCI$. Then a grouping of $\CFL_{in}(w)$ is obtained
by changing each $\PMCI$ with one of its groupings.

\begin{definition} \label{grouping1}
Let $\ell_1, \ldots , \ell_h$ be words in $L_{in}$ such that
$\ell_i$ is a prefix of $\ell_{i-1}$, $1 < i \leq h$.
We say that $(m_1, \ldots , m_k)$ is a {\rm grouping} of $(\ell_1, \ldots , \ell_h)$
if the following conditions are satisfied.
\begin{itemize}
\item[(1)]
$m_j$ is an inverse Lyndon word,
\item[(2)]
$\ell_1 \cdots  \ell_{h} = m_1 \cdots m_k$. More precisely, there are $i_0, i_1, \ldots , i_k$,
$i_0 = 0$, $1 \leq i_j \leq h$, $i_k = h$, such that $m_j = \ell_{i_{j-1} + 1} \cdots \ell_{i_j}$,
$1 \leq j \leq k$,
\item[(3)]
$m_1 \ll \ldots \ll m_k$.
\end{itemize}
\end{definition}

We now extend Definition \ref{grouping1} to $\CFL_{in}(w)$.

\begin{definition} \label{grouping2}
Let $w \in \Sigma^+$ and let $\CFL_{in}(w) = (\ell_1, \ldots , \ell_h)$.
We say that $(m_1, \ldots , m_k)$ is a {\rm grouping} of $\CFL_{in}(w)$ if
it can be obtained by replacing
any $\PMCI$ $\mathcal{C}$ in $\CFL_{in}(w)$
by a grouping of $\mathcal{C}$.
\end{definition}

Groupings of $\CFL_{in}(w)$
are inverse Lyndon factorizations of $w$ but there are
inverse Lyndon factorizations which are not groupings. Surprisingly
enough, $\ICFL(w)$ is a grouping of
$\CFL_{in}(w)$ but it is not always
its unique grouping.

\begin{proposition} \label{Pr3} \cite{inverseLyndon}
Let $(\Sigma, <)$ be a totally ordered alphabet. For any $w \in \Sigma^+$,
$\ICFL(w)$ is a grouping of $\CFL_{in}(w)$.
\end{proposition}

\begin{example} \label{groupings} \cite{inverseLyndon}.
{\rm Let $\Sigma = \{a,b,c,d \}$, $a < b < c < d$,
and $w = dabadabdabdadac$.
We have $\CFL_{in}(w) = (daba, dab, dab, dadac)$,
$\ICFL(w) = (daba, dabdab, dadac)$ (see Example \ref{nonunique}).
$\ICFL(w)$ is a grouping of $\CFL_{in}(w)$ but
$(dabadab, dabda, dac)$ is not a grouping.
Next, let $y = dabadabdabdabdadac$. We have
$\CFL_{in}(y) = (daba, dab, dab, dab, dadac)$
and $\ICFL(y) = (daba, (dab)^3, dadac)$.
The inverse Lyndon factorization $(dabadab, (dab)^2, dadac)$
is another grouping of $\CFL_{in}(y)$.}
\end{example}

Lemma \ref{bordo0} has been
proved in \cite{inverseLyndonTCS21}. It is needed in the proof
of Proposition \ref{bordo1}. In turn, Proposition \ref{bordo1} shows that the word $p$ in the pair
$(p, \bar{p}) \in \Pref_{bre}(w)$ has a grouping-like property.
Again, the proof of Proposition \ref{bordo1} was given in \cite{inverseLyndonTCS21} but
in \cite{inverseLyndonTCS21} uses a slightly incorrect statement of
Proposition \ref{preprime}.
However, the above-mentioned proof, also given below,
still holds because Proposition \ref{preprime} is in fact unnecessary
for it.

\begin{lemma} \label{bordo0}
Let $w \in \Sigma^+$ be a word which is not an inverse Lyndon word,
let $\CFL_{in}(w) = (\ell_1^{n_1}, \ldots , \ell_h^{n_h})$, with
$h  > 0$, $n_1, \ldots , n_h \geq 1$.
For all $z \in \Sigma^+$ and $b \in \Sigma$
such that $z$ is an anti-prenecklace, $zb$ is not an anti-prenecklace
and $zb$ is a prefix of $w$, there is an integer $g$ such that
$$zb = (u_1v_1)^{n_1} \cdots (u_gv_g)^{n_g} u_gb,$$
where $u_jv_j = u_ja_jv'_j = \ell_j, 1 \leq j \leq g, a_j < b$
and $u_gb$ is an anti-prenecklace.
\end{lemma}

\begin{remark} \label{twoborders} \cite{inverseLyndonTCS21}
Let $x, y$ two different borders of a same word $w \in \Sigma^+$.
If $x$ is shorter than $y$, then $x$ is a border of $y$.
\end{remark}

\begin{proposition} \label{bordo1}
Let $w \in \Sigma^+$ be a word which is not an inverse Lyndon word,
let $(p, \bar{p}) \in \Pref_{bre}(w)$ and let
$\ICFL(w) = (m_1, \ldots , m_k)$. Let
$\CFL_{in}(w) = (\ell_1^{n_1}, \ldots , \ell_h^{n_h})$, with
$h  > 0$, $n_1, \ldots , n_h \geq 1$ and let
$(\ell_1^{n_1}, \ldots , \ell_q^{n_q})$ be a
$\PMCI$ in $\CFL_{in}(w)$, $1 \leq q \leq h$.
Then the following properties hold.
\begin{itemize}
\item[(1)]
$p = \ell_1^{n_1} \cdots \ell_{g}^{n_{g}}$,
for some $g$, $1 \leq g \leq q$.
\item[(2)]
$\ell_{g} = u_g v_g = u_g a_g v'_g$, $\bar{p} = u_g b$, $a_g < b$.
\end{itemize}
\end{proposition}
\begdim
Let $w \in \Sigma^+$ be a word which is not an inverse Lyndon word,
let $(p, \bar{p}) \in \Pref_{bre}(w)$.
Let $\CFL_{in}(w) = (\ell_1^{n_1}, \ldots , \ell_h^{n_h})$, with
$h  > 0$, $n_1, \ldots , n_h \geq 1$ and let
$(\ell_1^{n_1}, \ldots , \ell_q^{n_q})$ be a
$\PMCI$ in $\CFL_{in}(w)$, $1 \leq q \leq h$.
Let $r, s \in \Sigma^*$, $a', b \in \Sigma$ be such that
$p = ra's$, $\bar{p} = rb$, $a' < b$.

By Proposition \ref{InvPreN}, the word $p \bar{p} = prb$ is not an anti-prenecklace but
its longest proper prefix is an anti-prenecklace.
Thus, by Lemma \ref{bordo0} there is an integer $g$ such that
$$p \bar{p} = (u_1v_1)^{n_1} \cdots (u_gv_g)^{n_g} u_gb,$$
where $u_jv_j = u_ja_jv'_j = \ell_j, 1 \leq j \leq g, a_j < b$
and $u_gb$ is an anti-prenecklace. Moreover $u_gb$ is a prefix of
$\ell_{g+1}^{n_{g+1}} \cdots  \ell_h^{n_h}$.
Let
$$\beta = (u_1v_1)^{n_1}(u_2v_2)^{n_2} \cdots (u_gv_g)^{n_g}.$$
The word $\beta$ is a nonempty proper prefix of $p \bar{p}$ thus, by
Definition \ref{brepref}, $\beta$ is an inverse Lyndon word. Therefore
$g \leq q$ (otherwise $g \geq q +1$, hence $\ell_q$ would be a prefix of $\beta$
and there would be a word $z'$ such that
$\ell_{q+1}z'$ is a suffix of $\beta$, a contradiction
since $\ell_q \ll \ell_{q+1}$ implies
$\beta \ll \ell_{q+1}z'$).
Moreover $\beta \ll u_gb$.

Since $p = ra's$, $\bar{p} = rb$, $a' < b$,
then $p \bar{p} = \beta u_gb = r a's rb$. By Proposition \ref{shortest}, $r$ is a suffix of
$u_g$. If $r = u_g$, then $p = \beta$ and the proof is ended. By contradiction, assume that
$r$ is a proper suffix of $u_g$.

Since $r$ is a proper suffix of $u_g$, we have
$|\beta| < |p| \leq |\beta u_g|$.
Hence, $\ell_g = u_gv_g$ and $u_g$ are nonempty prefixes of $p$.
By Lemma \ref{lem:inverse-Lyndon-word-prefix} $u_g$ is an inverse Lyndon word, thus,
by Proposition \ref{InvPreN}, $u_g$ is a nonempty anti-prenecklace.
Therefore, $u_g$ and $u_gb$ are both nonempty anti-prenecklaces.
By Theorem \ref{FundamentalPreneck},
there are $x, y \in \Sigma^*$, an integer $t \geq 1$, $c \in \Sigma$ such that
$xy$ is an anti-Lyndon word, $u_g = (xy)^t x$,
$y = cy'$ with $c \geq b$.

Observe that $m_1 = pr'$, for a prefix $r'$ of $r$. Thus
again $|\beta| < |m_1| \leq |\beta u_g|$. Moreover,
the words $\ell_{g+1}$ and $u_gb = (xy)^t xb$ are both prefixes of the
same word $\gamma = \ell_{g+1}^{n_{g+1}} \cdots  \ell_h^{n_h}$,
hence they are comparable for the prefix order. Since $\ell_{g+1}$
is the longest anti-Lyndon prefix of
$\gamma$, we have $|\ell_{g+1}| \geq |xy|$ and since $\ell_{g+1}$ is unbordered,
either $\ell_{g+1} = xy$ is a prefix of $\ell_{g}$ and $g+1 \leq q$,
or the word $u_gb = (xy)^t xb$ is a prefix of $\ell_{g+1}$.
By Proposition \ref{Pr3},
the first case holds, otherwise $m_1$ would not be a product
of anti-Lyndon words because $m_1$ is a prefix of $\beta u_g$
longer than $\beta$.

Recall that
$r$ is a proper suffix of $u_g$.
Moreover $r$ and $ra'$ are also prefixes of $u_g$
because $ra'$ and $u_g$ are both prefixes of $p$ and $r$ is shorter than
$u_g$. Therefore $r$ is a border of $u_g$ and $u_g$ starts with $ra'$.
Of course $r \not = x$ because $u_g$ starts with $ra'$ and also with $xc$,
with $c \geq b > a'$.
Now $r$ and $x$ are two different borders of $u_g$.
If $r$ were shorter than $x$, then $r$ would be a border of $x$
by Remark \ref{twoborders}.
This is impossible because $rcy'(xy)^{t-1}x$ would be a suffix of
the inverse Lyndon word $u_g$ and $u_g$ starts with $ra'$,
with $c \geq b > a'$. Thus $|r| > |x| \geq 0$.
Since $r$ is a nonempty border of $u_g = (xy)^t x$ and $|r| > |x| \geq 0$,
one of the following three cases holds:
\begin{eqnarray}
r & = & (xy)^{t'}x, \quad 0 < t' < t   \label{eq1} \\
r & = & y_1(xy)^{t'}x, \quad y_1 \mbox{ nonempty suffix of } y,  \quad 0 \leq  t' < t  \label{eq2} \\
r & = & x_1(yx)^{t'}, \quad x_1 \mbox{ nonempty suffix of } x,  \quad 0 < t' \leq t \label{eq3}
\end{eqnarray}
Assume that Eq. (\ref{eq1}) holds. Then $p$ starts with $ra' = (xy)^{t'}xa'$, $a' < b$,
and $p$ also starts with $u_g = (xy)^tx$. Since $t' < t$, the letter $a'$ should be the first
letter of $y = cy'$, $c \geq b > a'$. Therefore, Eq. (\ref{eq1}) cannot hold.

Assume that Eq. (\ref{eq2}) holds. Recall that $r$ is a prefix of $u_g$. Therefore $y_1$ is a prefix
of $xy$ and $y_1$ is also a nonempty suffix of $xy$. The word $xy$ is an anti-Lyndon word, thus $xy$ is umbordered.
Consequently, $y_1 = xy$, hence $x = 1$ and $y_1 =y$.
By Eq. (\ref{eq2}), we have $r = y^{t'+1}$, with $0 \leq  t' < t$. Moreover, $t'+1 < t$ since
$r$ is a proper suffix of $u_g = y^t$.
As above, $p$ starts with $ra' = y^{t'+1}a'$, $a' < b$,
and $p$ also starts with $u_g = y^t$. Since $t'+1 < t$, the letter $a'$ should be the first
letter of $y = cy'$, $c \geq b > a'$. Therefore, Eq. (\ref{eq2}) cannot hold.

Finally, assume that Eq. (\ref{eq3}) holds.
If $x_1 \not = x$, then $x_1y$ would be both a proper nonempty suffix
and a prefix of $xy$, hence a nonempty border of $xy$,
which is impossible since $xy$ is an anti-Lyndon word.
Therefore $x_1 = x$. If $t' < t$, then $r$ satisfies Eq. (\ref{eq1})
and we proved that this is impossible.
Thus $t' = t$, which implies $r = u_g$, a contradiction.
\enddim

The following result was proved in \cite{inverseLyndon} and will be used in Section
\ref{fromCtoI}.

\begin{proposition} \label{Pr0}
Let $(\Sigma, <)$ be a totally ordered alphabet. Let $w \in \Sigma^+$ and let
$\CFL_{in}(w) = (\ell_1, \ldots , \ell_h)$.
If $w$ is an inverse Lyndon word, then either $w$ is unbordered or
$\ell_1, \ldots , \ell_{h}$ is a $\PMCI$ in $\CFL_{in}(w)$.
In both cases $\ICFL(w) = (w)$ is the unique grouping of $\CFL_{in}(w)$.
\end{proposition}


\section{From $\ICFL$ to $\CFL_{in}$} \label{fromItoC}

In what follows, $(\Sigma, <)$ denotes a totally ordered alphabet.
Let $w \in \Sigma^+$ and let $\ICFL(w) = (m_1, \ldots , m_k)$, where
$m_i \ll m_{i+1}$, $1 \leq i \leq k-1$.
In this section we give an algorithm to construct $\CFL_{in}(w)$ starting from
$\ICFL(w)$ (Definition \ref{def:NB}).
To this aim, we first demonstrate that we can limit ourselves to handling the single factor
$m_i$ in $\ICFL(w)$ (Corollary \ref{cor2}).
This result is a direct consequence of the following
proposition.

\begin{proposition} \label{Pr4}
Let $w \in \Sigma^+$.
If $(m_1, \ldots , m_k)$ is a grouping of $\CFL_{in}(w)$,
then
$$\CFL_{in}(w) = (\CFL_{in}(m_1), \ldots , \CFL_{in}(m_k))$$
\end{proposition}
\begdim
Let $w \in \Sigma^+$ and let
$\CFL_{in}(w) = (\ell_1, \ldots , \ell_h)$.
We prove our statement by induction on $|w|$.
If $|w| = 1$, then there is only one grouping of $\CFL_{in}(w)$, namely
$(w)$. Obviously $\CFL_{in}(w) = (w)$ and we are done.

Assume $|w| > 1$.
Let $(m_1, \ldots , m_k)$ be a grouping of $\CFL_{in}(w)$
and let $w' \in \Sigma^*$ be such that $w = m_1 w'$.
According to Definition \ref{grouping2} there are words in $\CFL_{in}(w)$, say
$\ell_1, \ldots , \ell_v$, $1 \leq v \leq h$, such that
$\ell_i$ is a prefix of $\ell_{i-1}$, $1 < i \leq v$ and
$m_1 = \ell_1 \cdots \ell_v$.
Moreover, by Theorem \ref{Lyndon-factorization}
\begin{equation} \label{EQprima1}
\CFL_{in}(m_1) = (\ell_1, \ldots , \ell_v)
\end{equation}
If $w' = 1$, then $w = m_1$, $v = h$, hence
$\CFL_{in}(w) = \CFL_{in}(m_1)$ and we are done.

Otherwise, $w' = m_2 \cdots m_k \in \Sigma^+$, $v < h$ and,
by Theorem \ref{Lyndon-factorization},
$\CFL_{in}(w') = (\ell_{v+1}, \ldots , \ell_h)$.
In addition, $(m_2, \ldots , m_k)$ is a grouping of $\CFL_{in}(w')$
and $|w'| < |w|$. By using the induction hypothesis, we have
\begin{equation} \label{EQprima2}
\CFL_{in}(w') = (\ell_{v+1}, \ldots , \ell_h) = (\CFL_{in}(m_2), \ldots , \CFL_{in}(m_k))
\end{equation}
Finally, by Eqs. \eqref{EQprima1}, \eqref{EQprima2}, we obtain
\begin{equation}
\CFL_{in}(w) = (\ell_1, \ldots , \ell_h) = (\CFL_{in}(m_1), \ldots , \CFL_{in}(m_k))
\end{equation}
\enddim

\begin{corollary} \label{cor2}
Let $w \in \Sigma^+$.
If $\ICFL(w) = (m_1, \ldots , m_k)$,
then
$$\CFL_{in}(w) = (\CFL_{in}(m_1), \ldots , \CFL_{in}(m_k))$$
\end{corollary}
\begdim
Let $w \in \Sigma^+$ and let
$\ICFL(w) = (m_1, \ldots , m_k)$.
By Proposition \ref{Pr3}, $\ICFL(w)$ is a grouping of $\CFL_{in}(w)$, hence
our claim follows by Proposition \ref{Pr4}.
\enddim

In the following proposition, we show that each bordered word has a unique
nonempty border which is unbordered.

\begin{proposition} \label{Prdopo4}
Each bordered word has a unique
nonempty border which is unbordered.
\end{proposition}
\begdim
We first prove that each bordered word
$w \in \Sigma^+$ has a nonempty border $z$ which is unbordered,
by induction on the length of $w$.
Indeed, the shortest bordered word is $w = aa$,
where $a$ is a letter and $z = a$.

Let $|w| > 2$ and let $y$ be a nonempty border of $w$.
Thus, there are nonempty words $x, x'$
such that
\begin{equation} \label{EQprima3}
w = xy = y x'
\end{equation}
If $y$ is unbordered, we are done.
Otherwise, since $y$ is a bordered word shorter than
$w$, by induction hypothesis, $y$ has a
nonempty border $z$ which is unbordered and
there are $t, t' \in \Sigma^+$ such that
\begin{equation} \label{EQprima4}
y = tz = zt'
\end{equation}
By using Eqs. \eqref{EQprima3}, \eqref{EQprima4},
we get
\begin{equation} \notag
xtz = xy = w = yx' = zt'x'
\end{equation}
hence $z$ is a nonempty border of $w$ which is unbordered.

Now we prove that there is only one unbordered nonempty border of $w$.
Suppose that $y, v$ are two nonempty unbordered borders of $w$.
Thus there are $x, x', u, u' \in \Sigma^+$ such that
\begin{equation} \label{EQprima5}
w = xy = y x' = uv = v u'
\end{equation}
Assume $|y| \geq |v|$ (a similar argument holds if $|y| \leq |v|$).
By Eq. \eqref{EQprima5},
there are $t, t' \in \Sigma^*$ such that
\begin{equation} \notag
y = tv = v t'
\end{equation}
Since $y$ is unbordered we get $t = t' =1$, hence $y = v$.
\enddim

\begin{example}
{\rm Let $\Sigma = \{a, b \}$. The word $ababa$ has the nonempty borders $a, aba$
and $a$ is unbordered. The word $aaaa$ has the nonempty borders $a, aa, aaa$
and $a$ is unbordered.}
\end{example}

Definition \ref{def:NB} below provides a recursive construction of a sequence
$\NB(x)$ of suffixes of a word $x$.
Proposition \ref{Pr5} proves that every nonempty word $x$ is the concatenation of the elements of $\NB(x)$
and that there is only one sequence $\NB(x)$ associated with $x$ according to Definition \ref{def:NB}.
Subsequently Proposition \ref{Pr6} proves that if $w \in \Sigma^+$ and $x$ is a factor of a grouping of
of $\CFL_{in}(w)$, then $\CFL_{in}(x) = \NB(x)$.
According to Corollary \ref{cor2}, our aim, that is, to construct $\CFL_{in}(w)$ from $\ICFL_{in}(w)$,
is achieved (Corollary \ref{cor3}).

\begin{definition} \label{def:NB}
Let $x \in \Sigma^+$.
\begin{itemize}
\item []
{\rm (Basis Step)}
If $x$ is unbordered,
then $\NB(x) = (x)$.
\item []
{\rm (Recursive Step)}
If $x$ is bordered, let $z$ be the unique nonempty border of $x$ which is unbordered
and let $y \in \Sigma^+$ be such that $x = yz$.
Then, $\NB(x) = (\NB(y), z)$.
\end{itemize}
The {\rm length} $|\NB(x)|$ of $\NB(x)$ is the number of elements in $\NB(x)$.
\end{definition}

\begin{proposition} \label{Pr5}
For any word $x \in \Sigma^+$, there is a unique sequence
$\ell'_1, \ldots , \ell'_v$ of words over $\Sigma$
such that $\NB(x) = (\ell'_1, \ldots , \ell'_v)$.
If $\ell'_1, \ldots , \ell'_v$ is
such that $\NB(x) = (\ell'_1, \ldots , \ell'_v)$, then
$$x = \ell'_1 \cdots  \ell'_v$$
\end{proposition}
\begdim
Notice that if $x$ is unbordered, then, by Definition \ref{def:NB},
$\NB(x) = (x)$ is uniquely determined and we are done.
The proof is by induction on $|x|$.
By the above reasoning if $|x| = 1$, then the conclusion follows immediately
and the same holds for any unbordered word $x$.

Hence, let $x$ be a bordered word such that $|x| > 1$.
Let
\begin{equation} \label{EQprima6}
\NB(x) = (\ell'_1, \ldots , \ell'_v), \quad \NB(x) = (\mu_1, \ldots , \mu_p)
\end{equation}
with $v > 1$, $p > 1$.
Let $z$ be the unique nonempty border of $x$ which is unbordered
and let $y \in \Sigma^+$ be such that $x = yz$.
By Definition \ref{def:NB}, we have
\begin{equation} \label{EQprima8}
\NB(x)  =  (\NB(y), z)= (\ell'_1, \ldots , \ell'_v), \quad \NB(x)  =  (\NB(y), z) = (\mu_1, \ldots , \mu_p)
\end{equation}
Hence,
$\ell'_v = z  = \mu_p$ and $\NB(y) = (\ell'_1, \ldots , \ell'_{v - 1})$, $\NB(y) =
(\mu_1, \ldots , \mu_{p - 1})$.
Obviously $|y| < |x|$. By using induction hypothesis, we have
$p = v$, $\ell'_i = \mu_i$, $1 \leq i \leq v$,
$$x = yz = \ell'_1 \cdots  \ell'_{v-1} z = \ell'_1 \cdots  \ell'_v$$
and the proof is complete.
\enddim

The following example provides an interesting observation.

\begin{example}
{\rm Let $\Sigma = \{a,b,c,d \}$, $a < b < c < d$,
and $w = dabadabdabdadac$.
We have
$\ICFL(w) = (daba, dabdab, dadac)$ (see Example \ref{nonunique}).
Then, $\NB(daba) = daba$, $\NB(dabdab) = (dab , dab)$,
$\NB(dadac) = dadac$.
Notice that
\begin{eqnarray*}
\CFL_{in}(w) & = & (daba, dab, dab, dadac) = (\NB(daba), \NB(dabdab), \NB(dadac)) \\
             & \not = & (dabadabdabdadac) = (\NB(dabadabdabdadac)) = (\NB(w))
\end{eqnarray*}
Next, let $y = dabadabdabdabdadac$. We have
$\ICFL(w) = (daba, (dab)^3, dadac)$
and
$$\CFL_{in}(y) = (daba, dab, dab, dab, dadac) =
(\NB(daba), \NB((dab)^3), \NB(dadac))$$}
\end{example}

\begin{proposition} \label{Pr6}
Let $(\Sigma, <)$ be a totally ordered alphabet.
Let $m, \ell_1, \ldots , \ell_h$ be words in $\Sigma^+$
such that
\begin{itemize}
\item[(1)]
$m$ is an inverse Lyndon word.
\item[(2)]
$\ell_1, \ldots , \ell_h$ are words in $L_{in}$ such that
$\ell_i$ is a prefix of $\ell_{i-1}$, $1 < i \leq h$.
\item[(3)]
$m =\ell_1 \cdots  \ell_{h}$.
\end{itemize}
Then the two sequences $\CFL_{in}(m)$, $\NB(m)$ are equal, that is,
$\CFL_{in}(m) = \NB(m)$.
\end{proposition}
\begdim
The proof is by induction on $|\NB(m)|$.
If $|\NB(m)| = 1$, then by Definition \ref{def:NB},
the inverse Lyndon word $m$ is unbordered and $\NB(m) = (m)$.
Thus, by Proposition \ref{Pr2}, the word $m$ is
is an anti-Lyndon word.
Consequently, $\CFL_{in}(m) = (m) = \NB(m)$.

Otherwise, let $|\NB(m)| > 1$. Let $y, z$ be such that
$z$ is the unique unbordered border of $m$ and $m = yz$.
Then, by Definition \ref{def:NB}, $\NB(m) = (\NB(y), z)$.
Let $\ell_1, \ldots , \ell_h$ be as in the statement.
Of course $\CFL_{in}(m) = (\ell_1, \ldots , \ell_h)$.
Let $\NB(y) = (\ell'_1, \ldots  , \ell'_v)$.
By item (3) and by Proposition \ref{Pr5} applied to $y$, we have
\begin{equation} \label{EQ1}
m = yz = \ell'_1 \cdots  \ell'_v z = \ell_1 \cdots  \ell_h
\end{equation}
Now, $\ell_h$ is a border of $m$, because $\ell_h$ is a prefix
of $\ell_1$, hence of $m$. Moreover, by Proposition \ref{Pr1},
$\ell_h$ is an unbordered word.
By Proposition \ref{Prdopo4} and looking at Eq. \eqref{EQ1}, we have
\begin{equation} \label{EQ2}
z = \ell_h, \quad y = \ell'_1 \cdots  \ell'_v = \ell_1 \cdots  \ell_{h - 1}
\end{equation}
The word $y$ is an inverse Lyndon word because it is a nonempty prefix of the inverse Lyndon
word $m$ (Lemma \ref{lem:inverse-Lyndon-word-prefix}). Moreover
$\CFL_{in}(y) = (\ell_1, \ldots , \ell_{h - 1})$ and
$y, \ell_1, \ldots  , \ell_{h- 1}$ satisfy
the hypothesis of the statement with $|\NB(y)| < |\NB(m)|$.
Thus, by using the induction hypothesis and Eq. \eqref{EQ2},
$\CFL_{in}(y)$ and $\NB(y)$ are equal, that is,
\begin{equation} \label{EQ3}
(\ell_1, \ldots , \ell_{h - 1}) = \CFL_{in}(y) = \NB(y) = (\ell'_1, \ldots , \ell'_v)
\end{equation}
and consequently
\begin{equation} \label{EQ4}
h-1 = v, \quad \ell_i = \ell'_i, \quad 1 \leq i \leq h-1
\end{equation}
Finally, by Eqs. \eqref{EQ2}, \eqref{EQ3}, \eqref{EQ4}, we have
$$\CFL_{in}(m) = (\ell_1, \ldots , \ell_h) = (\ell'_1, \ldots , \ell'_v, \ell_h) =
(\ell'_1, \ldots , \ell'_v , z) = (\NB(y), z) = \NB(m)$$
\enddim

\begin{proposition} \label{Pr7}
Let $w \in \Sigma^+$.
If $(m_1, \ldots , m_k)$ is a grouping of $\CFL_{in}(w)$,
then for each $i$, $1 \leq i \leq v$,
$\CFL_{in}(m_i) = \NB(m_i)$.
\end{proposition}
\begdim
Let $w \in \Sigma^+$, let
$(m_1, \ldots , m_k)$ be a grouping of $\CFL_{in}(w)$.
Looking at Definitions \ref{grouping1}, \ref{grouping2},
the conclusion follows easily from
Proposition \ref{Pr6} applied to the word $m_i$,
for each $i$, $1 \leq i \leq v$.
\enddim

\begin{proposition} \label{Pr8}
Let $w \in \Sigma^+$.
If $(m_1, \ldots , m_k)$ is a grouping of $\CFL_{in}(w)$,
then $\CFL_{in}(w) = (\NB(m_1), \ldots , \NB(m_k))$.
\end{proposition}
\begdim
Let $w \in \Sigma^+$, let
$(m_1, \ldots , m_k)$ be a grouping of $\CFL_{in}(w)$.
By Propositions \ref{Pr4}, \ref{Pr7}, we have
$$\CFL_{in}(w) = (\CFL_{in}(m_1), \ldots , \CFL_{in}(m_k)) = (\NB(m_1), \ldots , \NB(m_k))$$
\enddim

\begin{corollary} \label{cor3}
Let $w \in \Sigma^+$.
If $\ICFL(w) = (m_1, \ldots , m_k)$,
then
$$\CFL_{in}(w) = (\NB(m_1), \ldots , \NB(m_k))$$
\end{corollary}
\begdim
The conclusion follows easily from
Propositions \ref{Pr3} and \ref{Pr8}.
\enddim

\begin{example}
{\rm Let $\Sigma = \{a,b,c,d \}$, $a < b < c < d$,
and let $y = dabadabdabdabdadac$.
We know that $\ICFL(y) = (daba, (dab)^3, dadac)$
(see Example \ref{groupings}). By Corollaries \ref{cor2} and \ref{cor3},
we have
\begin{eqnarray*}
\CFL_{in}(y) & = & (\CFL_{in}(daba), \CFL_{in}((dab)^3), \CFL_{in}(dadac)) \\
             & = & (\NB(daba), \NB((dab)^3), \NB(dadac)) \\
             & = & (daba, dab, dab, dab, dadac)
\end{eqnarray*}
The inverse Lyndon factorization $(dabadab, (dab)^2, dadac)$
is another grouping of $\CFL_{in}(y)$.
By Propositions \ref{Pr4} and \ref{Pr8}, we have
\begin{eqnarray*}
\CFL_{in}(y) & = & (\CFL_{in}(dabadab), \CFL_{in}((dab)^2), \CFL_{in}(dadac)) \\
               & = & (\NB(dabadab), \NB((dab)^2), \NB(dadac)) \\
               & = & (daba, dab, dab, dab, dadac)
\end{eqnarray*}
}
\end{example}


\section{From $\CFL_{in}$ to $\ICFL$} \label{fromCtoI}

In what follows, $(\Sigma, <)$ denotes a totally ordered alphabet.
Let $w \in \Sigma^+$.
The aim of this section is to show how to get $\ICFL(w)$ from $\CFL_{in}(w)$.
As we know,
$\ICFL(w) = (m_1, \ldots , m_k)$ is a grouping of
$\CFL_{in}(w)$ (Proposition \ref{Pr3}).
We also know that the word $p$ in the pair
$(p, \bar{p}) \in \Pref_{bre}(w)$ has a grouping-like property
(Proposition \ref{bordo1}).
These properties allow us to prove that $\ICFL(w)$
can be computed locally in the non-increasing maximal chains for the prefix order
of its decomposition (Proposition \ref{step1}).

\begin{remark}
{\rm Let $(\ell'_1, \ldots , \ell'_v)$ be a non-increasing chain for the prefix order
of anti-Lyndon words. Let $y = \ell'_1 \cdots  \ell'_v$.
It is easy to see that $\CFL_{in}(y) = (\ell'_1, \ldots , \ell'_v)$.}
\end{remark}

\begin{proposition} \label{step1}
Let $w \in \Sigma^+$,
let $\CFL_{in}(w) = (\ell_1, \ldots , \ell_h)$ and
let $(\mathcal{C}_{1}, \mathcal{C}_{2}, \ldots , \mathcal{C}_{s})$ be
the decomposition of $\CFL_{in}(w)$ into its non-increasing maximal chains for the prefix order.
Let $w_1, \ldots , w_s$ be words such that
$\CFL_{in}(w_j) = \mathcal{C}_{j}$, $1 \leq j \leq s$.
Then $\ICFL(w)$ is the concatenation of the sequences
$\ICFL(w_1), \ldots , \ICFL(w_s)$, that is,
\begin{equation} \label{EqLocal}
\ICFL(w) = (\ICFL(w_1), \ldots , \ICFL(w_s))
\end{equation}
\end{proposition}
\begdim
Let $w \in \Sigma^+$,
let $\CFL_{in}(w) = (\ell_1, \ldots , \ell_h)$ and
let $(\mathcal{C}_{1}, \mathcal{C}_{2}, \ldots , \mathcal{C}_{s})$ be
the decomposition of $\CFL_{in}(w)$ into its non-increasing maximal chains for the prefix order.
Let $w_1, \ldots , w_s$ be words such that
$\CFL_{in}(w_j) = \mathcal{C}_{j}$, $1 \leq j \leq s$.
Let $\ICFL(w) = (m_1, \ldots , m_k)$.
The proof is by induction on $|w|$.
If $|w| = 1$, then $w$ is an inverse Lyndon word
and by Proposition \ref{Pr0} we are done.
Therefore assume $|w| > 1$.
If $k = 1$, by Definition \ref{def:ICFL}, $w = m_1$ is an inverse Lyndon word
and by Proposition \ref{Pr0} we are done.
Assume $k > 1$, thus $w$ is a word which is not an inverse Lyndon word.

By Proposition \ref{Pr3}, there are indexes $i_2, \ldots , i_s$,
$1 < i_2   \ldots  < i_s \leq k$
such that
$$w_1 = m_1 \cdots m_{i_2 -1}, \; w_2 = m_{i_2} \cdots m_{i_3 -1}, \ldots , w_s = m_{i_s} \cdots m_k$$
We actually prove more, namely we prove Eq. \eqref{EqLocalBis} below:
\begin{equation} \label{EqLocalBis}
\ICFL(w_1) = (m_1, \ldots , m_{i_2 -1}), \; \ICFL(w_2) = (m_{i_2}, \ldots,  m_{i_3 -1}),  \ldots , \ICFL(w_s) = (m_{i_s}, \ldots , m_k)
\end{equation}

Notice that if $w_1$ is an inverse Lyndon word, then $w_1 = m_1$. Otherwise, since $m_1 \ll m_2$ and
by item (2) in Lemma \ref{proplexord}, we would have
$w_1 = m_1 \cdots m_{i_2 -1} \ll m_2 \cdots m_{i_2 -1}$, a contradiction.
By a similar argument, if $w_1$ is an inverse Lyndon word, then
$w_1w_2$ is not an inverse
Lyndon word. Indeed, let $\ell$ the last element in $\mathcal{C}_{1}$
and $\ell'$ the first element in $\mathcal{C}_{2}$. The word $\ell$ is a prefix
of $w_1$ and by $\ell \ll \ell'$ we get $w_1w_2 \ll w_2$.

Let $(p,\overline{p})$ be the canonical pair associated with $w$ and let
$v \in \Sigma^*$ be such that $w = pv$. If $w_1$ is not an inverse Lyndon word,
let $(q,\overline{q})$ be the canonical pair associated with $w_1$.
The word $q\overline{q}$ is a prefix of $w$ which is not an inverse Lyndon word.
Since $p\overline{p}$ is the shortest prefix of $w$ which is not an inverse Lyndon word,
we have $|p\overline{p}| \leq |q\overline{q}|$,
hence $p\overline{p}$ is a prefix of $w_1$.
In turn, since $q\overline{q}$ is the shortest prefix of $w_1$ which is not an inverse Lyndon word,
we have $|q\overline{q}| \leq |p\overline{p}|$.
In conclusion, $p\overline{p} = q\overline{q}$ and, looking at Proposition \ref{shortest},
we see that $p = q$ and $\overline{p} = \overline{q}$. If $w_1$ is an inverse Lyndon word,
the same reasoning applies to
the canonical pair $(q,\overline{q})$ associated with $w_1w_2$.

Let $\ICFL(v) = (m'_1, \ldots, m'_{k})$ and let
$r,s \in \Sigma^*$, $a,b \in \Sigma$ such that $p = ras$, $\overline{p} = rb$ with $a < b$.
By Definition \ref{def:ICFL} one has
\begin{equation} \label{eqICFL}
\ICFL(w) = \begin{cases} (p, \ICFL(v)) & \mbox{ if } \overline{p} = rb \leq_{p} m'_1 \\
(p m'_1, m'_2, \ldots, m'_{k}) & \mbox{ if } m'_1 \leq_{p} r \end{cases}
\end{equation}
By Proposition \ref{bordo1}, there is $g$, $1 \leq g \leq h$, such that
$p = \ell_1 \cdots  \ell_g$. Let $w'_1 \in \Sigma^*$ be
such that $w_1 = p w'_1$. If $w'_1 \not = 1$, let $\mathcal{C}'_{1} = \CFL_{in}(w'_1)$.
Since $\mathcal{C}'_{1}$ is $\mathcal{C}_{1}$ after erasing $\ell_1, \ldots , \ell_g$, it
is easy to see that $\mathcal{C}'_{1}$ is a non-increasing maximal chain for the prefix order.
Hence the decomposition $\mathcal{D}$ of $\CFL_{in}(v)$ into its non-increasing maximal chains for the prefix order
is
$$\mathcal{D} = \begin{cases} (\mathcal{C}'_{1}, \mathcal{C}_{2}, \ldots , \mathcal{C}_{s}) & \mbox{ if } w'_1 \not = 1 \\
                              (\mathcal{C}_{2}, \ldots , \mathcal{C}_{s})& \mbox{ if } w'_1 = 1 \end{cases}$$
Of course $|v| < |w|$, thus by induction hypothesis, we have
\begin{equation} \label{EqLocalTris}
\ICFL(v) = \begin{cases} (\ICFL(w'_1), \ldots , \ICFL(w_s)) & \mbox{ if } w'_1 \not = 1 \\
                         (\ICFL(w_2), \ldots , \ICFL(w_s)) & \mbox{ if } w'_1 = 1 \end{cases}
\end{equation}
More specifically,
\begin{equation} \label{EqLocal4}
\ICFL(w_2) = (m_{i_2}, \ldots,  m_{i_3 -1}),  \ldots , \ICFL(w_s) = (m_{i_s}, \ldots , m_k)
\end{equation}
Moreover, if $w'_1 \not = 1$ one has
\begin{equation} \label{EqLocal5}
\ICFL(w'_1) = \begin{cases} (m_2, \ldots, m_{i_2 -1}) & \mbox{ if } m_1 = p \\
                         (m'_1, m_2, \ldots, m_{i_2 -1}) & \mbox{ if } m_1 = p m'_1 \end{cases}
\end{equation}
Notice that if $m_1 = p$ and $w'_1 \not = 1$, then $w_1$ is not an inverse Lyndon word (because
$w_1 \not = m_1$) and, by Definition \ref{def:ICFL},
$(p, \ICFL(w'_1)) = \ICFL(w_1)$ (because the canonical pair associated with $w_1$ is equal to the canonical
pair associated with $w$).
Furthermore, by Eqs. \eqref{eqICFL}, \eqref{EqLocalTris}, \eqref{EqLocal4},
$\ICFL(w) = (p, \ICFL(v)) = (p, \ICFL(w'_1), \ldots , \ICFL(w_s)) = (\ICFL(w_1), \ldots , \ICFL(w_s))$
Thus, Eqs. \eqref{EqLocal}, \eqref{EqLocalBis} hold.
Analogously, by Definition \ref{def:ICFL}, if $m_1 = p$ and $w'_1 = 1$, then
$w_1 = p$ and clearly $\ICFL(w_1) = (p) = (m_1)$. Thus, by Eqs. \eqref{EqLocalTris}, \eqref{EqLocal4},
$\ICFL(w) = (p, \ICFL(v)) = (p, \ICFL(w_2), \ldots , \ICFL(w_s)) = (\ICFL(w_1), \ICFL(w_2), \ldots , \ICFL(w_s))$
and Eqs. \eqref{EqLocal}, \eqref{EqLocalBis} hold.

Otherwise, if $m_1 = p m'_1$, then $w'_1 \not = 1$.
By Definition \ref{def:ICFL} and Eq. \eqref{EqLocal5},
we have $\ICFL(w_1) = (pm'_1, m_2, \ldots m_{i_2 -1}) = (m_1, m_2, \ldots m_{i_2 -1})$.
Thus, by Eqs. \eqref{eqICFL}, \eqref{EqLocalTris}, \eqref{EqLocal4}, \eqref{EqLocal5}, we have
\begin{eqnarray*}
\ICFL(w) & = & (m_1, m_2, \ldots , m_k) \\
         & = & (pm'_1, m_2, \ldots m_{i_2 -1}, m_{i_2} \ldots m_{i_3 -1}, \ldots , m_{i_s} \cdots m_k \\
         & = & (pm'_1, m_2, \ldots m_{i_2 -1}, \ICFL(w_2), \ldots , \ICFL(w_s)) \\
         & = & (\ICFL(w_1), \ICFL(w_2), \ldots , \ICFL(w_s))
\end{eqnarray*}
and Eqs. \eqref{EqLocal}, \eqref{EqLocalBis} hold.
This ends the proof.
\enddim

Proposition \ref{step1} shows that to obtain $\ICFL(w)$ from $\CFL_{in}(w)$ we can limit ourselves to the case
in which $\CFL_{in}(w)$ is a chain with respect to the prefix order.
Thus in the results that follow we will focus on these chains.
In Proposition \ref{InW1} we will prove that some products of consecutive elements
in such a chain form an inverse Lyndon word.
Then we will prove some properties of such products in Propositions \ref{InW2},
\ref{InW3}. We will use these properties to establish which of the aforementioned
products can be elements of an inverse Lyndon factorization (Proposition \ref{InF1}).

\begin{proposition} \label{InW1}
Let $\ell_1, \ldots , \ell_h$ be anti-Lyndon words that form
a non-increasing chain for the prefix order, that is,
$$\ell_{1} \geq_p \ell_{2} \geq_p \ldots  \geq_p \ell_{h}$$
Let $h \geq 2$ and let $i$, $j$, $1 \leq i < j \leq h$, be such that
\begin{equation} \notag
\ell_{1} = \ell_{2} = \ldots =  \ell_i \not = \ell_{i + 1}
\end{equation}
and $\ell_{i + 1} \cdots  \ell_j$ is a prefix of $\ell_1$.
Then, the word $\ell_1 \ell_2 \cdots \ell_{i +1} \cdots  \ell_j$
is an inverse Lyndon word.
\end{proposition}
\begdim
Let $\ell_1, \ldots , \ell_h$ be anti-Lyndon words that form
a non-increasing chain for the prefix order, that is,
$$\ell_{1} \geq_p \ell_{2} \geq_p \ldots  \geq_p \ell_{h}$$

Let $h \geq 2$ and let $i$, $j$, $1 \leq i < j \leq h$, be such that
\begin{equation} \notag
\ell_{1} = \ell_{2} = \ldots =  \ell_i \not = \ell_{i + 1}
\end{equation}
and $\ell_{i + 1} \cdots  \ell_j$ is a prefix of $\ell_1$.
Let $u, v$ be words such that
$$u = \ell_{i + 1} \cdots  \ell_j, \quad \ell_1 = uv $$
Thus
$$\ell_1 \ell_2 \cdots \ell_{i + 1} \cdots  \ell_j = (uv)^iu$$
Therefore the word
$\ell_1 \ell_2 \cdots \ell_{i + 1} \cdots  \ell_j$
is a sesquipower of the anti-Lyndon word $\ell_1$,
hence, by Proposition \ref{InvPreN},
$\ell_1 \ell_2 \cdots \ell_{i + 1} \cdots  \ell_j$ is an inverse Lyndon word.
\enddim

\begin{proposition} \label{InW2}
Let $\ell_1, \ldots , \ell_h$ be anti-Lyndon words that form
a non-increasing chain for the prefix order, that is,
$$\ell_{1} \geq_p \ell_{2} \geq_p \ldots  \geq_p \ell_{h}$$
Let $h \geq 2$ and let $i$, $1 \leq i < h$, be such that
\begin{equation} \notag
\ell_{1} = \ell_{2} = \ldots =  \ell_i \not = \ell_{i + 1}
\end{equation}
If there is $j$, $i < j \leq h$, such that $\ell_{i + 1} \cdots  \ell_j$ is a prefix of $\ell_1$,
then $|\ell_1| > |\ell_{i + 1} \cdots  \ell_j|$.
\end{proposition}
\begdim
Let $\ell_1, \ldots , \ell_h$ be anti-Lyndon words that form
a non-increasing chain for the prefix order, that is,
$$\ell_{1} \geq_p \ell_{2} \geq_p \ldots  \geq_p \ell_{h}$$
Let $h \geq 2$ and let $i$, $j$, $1 \leq i < j \leq h$, be such that
\begin{equation} \notag
\ell_{1} = \ell_{2} = \ldots =  \ell_i \not = \ell_{i + 1}
\end{equation}
and $\ell_{i + 1} \cdots  \ell_j$ is a prefix of $\ell_1$.
Thus $|\ell_1| \geq |\ell_{i +1} \cdots  \ell_j|$.
By contradiction, assume $|\ell_1| = |\ell_{i +1} \cdots  \ell_j|$.
Hence $\ell_1 = \ell_{i +1} \cdots  \ell_j$.
In addition $|\ell_t| < |\ell_1|$, $i + 1 \leq t \leq j$,
because $\ell_t$ is a prefix of $\ell_{t-1}$
and $\ell_1 = \ell_i \not = \ell_{i + 1}$.
Consequently, $\ell_j$ would be a nonempty border of $\ell_1$,
in contradiction with Proposition \ref{Pr1}.
\enddim

\begin{proposition} \label{InW3}
Let $\ell_1, \ldots , \ell_h$ be anti-Lyndon words that form
a non-increasing chain for the prefix order, that is,
$$\ell_{1} \geq_p \ell_{2} \geq_p \ldots  \geq_p \ell_{h}$$
Let $h \geq 2$ and let $i$, $1 \leq i < h$, be such that
\begin{equation} \notag
\ell_{1} = \ell_{2} = \ldots =  \ell_i \not = \ell_{i + 1}
\end{equation}
If there is $j$, $i < j \leq h$, such that $\ell_{i + 1} \cdots  \ell_j$ is a prefix of
$\ell_1 \ell_2 \cdots \ell_{i}$,
then $\ell_{i + 1} \cdots  \ell_j$ is a prefix of $\ell_1$.
\end{proposition}
\begdim
Let $\ell_1, \ldots , \ell_h$ be anti-Lyndon words that form
a non-increasing chain for the prefix order, that is,
$$\ell_{1} \geq_p \ell_{2} \geq_p \ldots  \geq_p \ell_{h}$$
Let $h \geq 2$ and let $i$, $j$, $1 \leq i < j \leq h$, be such that
\begin{equation} \notag
\ell_{1} = \ell_{2} = \ldots =  \ell_i \not = \ell_{i + 1}
\end{equation}
and $\ell_{i + 1} \cdots  \ell_j$ is a prefix of $\ell_1 \ell_2 \cdots \ell_{i}$.
If $\ell_{i + 1} \cdots  \ell_j$ were not a prefix of $\ell_1$, then $i > 1$, $j > i + 1$
and there would exist $q$, $i + 1 < q \leq j$, such that
$$\ell_1 = \ell_{i +1} \cdots  \ell'_q, \quad  \ell_q = \ell'_q \ell''_q, \quad \ell'_q \not = 1$$
If $\ell''_q = 1$, then $\ell_q = \ell'_q$ would be a proper suffix of $\ell_1$
which is also a nonempty prefix of $\ell_1$, hence $\ell_q$ would be a nonempty border of $\ell_1$,
in contradiction with Proposition \ref{Pr1}. Consequently $\ell''_q \not = 1$ would be a proper suffix of $\ell_q$
which is also a nonempty prefix of $\ell_2$ and therefore of $\ell_q$.
Hence $\ell''_q$ would be a nonempty border of $\ell_q$,
in contradiction with Proposition \ref{Pr1}.
\enddim

\begin{proposition} \label{InF1}
Let $\ell_1, \ldots , \ell_h$ be anti-Lyndon words that form
a non-increasing chain for the prefix order, that is,
$$\ell_{1} \geq_p \ell_{2} \geq_p \ldots  \geq_p \ell_{h}$$
Let $h \geq 2$ and let $i$, $j$, $1 \leq i < j < h$, be such that
\begin{equation} \notag
\ell_{1} = \ell_{2} = \ldots =  \ell_i \not = \ell_{i + 1}
\end{equation}
and $\ell_{i + 1} \cdots  \ell_j$ is a prefix of $\ell_1$.
If $\ell_{i + 1} \cdots  \ell_{j +1}$ is not
a prefix of $\ell_{1}$, then
one has
\begin{equation} \notag
\ell_{1} \ll \ell_{i + 1} \cdots  \ell_{j +1}
\end{equation}
More specifically, there are words $r, s, s' \in \Sigma^*$ and $a, b \in \Sigma$, $a < b$ such that
\begin{equation} \notag
\ell_{1} = \ell_{i + 1} \cdots  \ell_j ras, \quad
\ell_{j + 1} = rbs'
\end{equation}
\end{proposition}
\begdim
Let $\ell_1, \ldots , \ell_h$ be anti-Lyndon words that form
a non-increasing chain for the prefix order, that is,
$$\ell_{1} \geq_p \ell_{2} \geq_p \ldots  \geq_p \ell_{h}$$
Let $h \geq 2$ and let $i$, $j$, $1 \leq i < j < h$, be such that
\begin{equation} \notag
\ell_{1} = \ell_{2} = \ldots =  \ell_i \not = \ell_{i + 1}
\end{equation}
and $\ell_{i + 1} \cdots  \ell_j$ is a prefix of $\ell_1$.

By Proposition \ref{InW2} one has $|\ell_1| > |\ell_{i + 1} \cdots  \ell_j|$,
hence there is a nonempty word $x$ such that $\ell_1 = \ell_{i + 1} \cdots  \ell_j x$.
Notice that $x$ cannot be a prefix of $\ell_{j + 1}$ because $\ell_{j + 1}$
is a proper prefix of $\ell_1$ and $\ell_1$ is unbordered (Proposition \ref{Pr1}).
If $\ell_{i+1} \cdots  \ell_{j +1 }$ is not
a prefix of $\ell_{1}$, then $\ell_{j + 1}$ cannot be a prefix of $x$ either.

Thus, there are words $r, s, s' \in \Sigma^*$ and $a, b \in \Sigma$, $a \not = b$ such that
\begin{equation} \notag
\ell_{1} = \ell_{i + 1} \cdots  \ell_j ras, \quad
\ell_{j + 1} = rbs'
\end{equation}
The word $\ell_{j + 1}$ is a prefix of $\ell_{i + 1}$, hence $\ell_{j + 1}$ is a prefix
of $\ell_{1}$ such that $|\ell_{j + 1}| \leq |\ell_{i + 1}|$.
Therefore, there is $z \in \Sigma^*$
such that
$\ell_1 = rbs'zras$ which yields $a < b$ because
$\ell_1$ is an inverse Lyndon word.
Consequently,
$$\ell_{1} = \ell_{i + 1} \cdots  \ell_j ras \ll
\ell_{i + 1} \cdots  \ell_j rbs' = \ell_{i + 1} \cdots \ell_{j + 1}$$
\enddim

In Propositions \ref{Algo1} - \ref{Algo3} we prove some properties of the
factors in $\CFL_{in}$.
At the end of this section we briefly discuss how
Propositions \ref{Algo2}, \ref{Algo3}
can be used to obtain  $\ICFL$ from
$\CFL_{in}$.

\begin{proposition} \label{Algo1}
Let $w \in \Sigma^+$, let $(\ell_1, \ldots , \ell_h)$
be a chain for the prefix order
in $\CFL_{in}(w)$.
The word $\ell_{1} \cdots  \ell_h$ is not an inverse Lyndon word if and only
if $h > 2$ and there are $i$, $j$, $1 \leq i < j < h$, such that
\begin{itemize}
\item[(1)]
$\ell_{1} = \ell_{2} = \ldots =  \ell_i \not = \ell_{i + 1}$,
\item[(2)]
$\ell_{i + 1} \cdots  \ell_j$ is a prefix of $\ell_1$
\item[(3)]
$\ell_{i + 1} \cdots  \ell_{j +1}$ is not
a prefix of $\ell_{1}$
\item[(4)]
We have
\begin{equation} \notag
\ell_{1} \ll \ell_{i + 1} \cdots  \ell_{j +1}
\end{equation}
More specifically, there are words $r, s, s' \in \Sigma^*$ and $a, b \in \Sigma$, $a < b$ such that
\begin{equation} \notag
\ell_{1} = \ell_{i + 1} \cdots  \ell_j ras, \quad
\ell_{j + 1} = rbs'
\end{equation}
\end{itemize}
\end{proposition}
\begdim
Let $w \in \Sigma^+$, let $(\ell_1, \ldots , \ell_h)$
be a chain for the prefix order
in $\CFL_{in}(w)$.
Let $h > 2$ and let $i$, $j$, $1 \leq i < j < h$, be such that
items (1)-(4) in the statement are satisfied. We prove that
$\ell_{1} \cdots  \ell_h$ is not an inverse Lyndon word.
Indeed, by item (4)
$$\ell_{1} \ll \ell_{i + 1} \cdots  \ell_{j +1}$$
Thus, by item (2) in Lemma \ref{proplexord},
$$\ell_{1} \ell_2 \cdots  \ell_h \ll \ell_{i + 1} \cdots  \ell_{j +1} \ell_{j + 2} \cdots  \ell_h$$
therefore $\ell_{1} \ell_2 \cdots  \ell_h$ does not satisfy Definition \ref{inverse-Lyndon-word}.

Conversely, if $\ell_{1} \cdots  \ell_h$ is not an inverse Lyndon word, then
there is $i$, $1 \leq i < h$, such that
$$\ell_{1} = \ell_{2} = \ldots =  \ell_i \not = \ell_{i + 1},$$
otherwise $\ell_{1} \cdots  \ell_h = (\ell_1)^h$ would be a sesquipower of the
anti-Lyndon word $\ell_1$, hence, by Proposition \ref{InvPreN},
$\ell_{1} \cdots  \ell_h$ would be an inverse Lyndon word.
Furthermore, there is $j$, $i < j < h$, such that
$\ell_{i + 1} \cdots  \ell_j$ is a prefix of $\ell_1$ but
$\ell_{i + 1} \cdots  \ell_{j + 1}$ is not
a prefix of $\ell_{1}$ since otherwise
$\ell_{i + 1} \cdots  \ell_h$ would be a prefix of $\ell_1$
and $\ell_{1} \cdots  \ell_h$ would be an inverse Lyndon word
by Proposition \ref{InW1}.
Thus, by Proposition \ref{InF1},
\begin{equation} \notag
\ell_{1} \ll \ell_{i + 1} \cdots  \ell_{j +1}
\end{equation}
More specifically, there are words $r, s, s' \in \Sigma^*$ and $a, b \in \Sigma$, $a < b$ such that
\begin{equation} \notag
\ell_{1} = \ell_{i + 1} \cdots  \ell_j ras, \quad
\ell_{j + 1} = rbs'
\end{equation}
and the proof is complete.
\enddim

\begin{proposition} \label{Algo2}
Let $w \in \Sigma^+$, let $(\ell_1, \ldots , \ell_h)$
be a non-increasing maximal chain for the prefix order ($\PMCI$)
in $\CFL_{in}(w)$, where $\ell_1$ is the first element in
$\CFL_{in}(w)$, and let $\ICFL(w) = (m_1, \ldots , m_k)$.
If $\ell_{1} \cdots  \ell_h$ is an inverse Lyndon word, then
$$m_1 = \ell_{1} \cdots  \ell_h$$
\end{proposition}
\begdim
Let $w \in \Sigma^+$, let $(\ell_1, \ldots , \ell_h)$
be a non-increasing maximal chain for the prefix order ($\PMCI$)
in $\CFL_{in}(w)$, where $\ell_1$ is the first element in
$\CFL_{in}(w)$.
If $\ell_{1} \cdots  \ell_h$ is an inverse Lyndon word, then
$m_1 = \ell_{1} \cdots  \ell_h$. Otherwise, by Proposition \ref{Pr3},
there would exist $j$, $t$, $1 \leq j < t \leq h$, such that
$$\ell_{1} \cdots  \ell_j = m_1 \ll m_2 = \ell_{j + 1} \cdots  \ell_t$$
and this would imply, by item (2) in Lemma \ref{proplexord},
$$\ell_{1} \cdots  \ell_h = (\ell_{1} \cdots  \ell_j) \ell_{j + 1} \cdots  \ell_h
\ll \ell_{j + 1} \cdots  \ell_h$$
in contradiction with Definition \ref{inverse-Lyndon-word}.
\enddim

\begin{proposition} \label{Algo3}
Let $w \in \Sigma^+$, let $(\ell_1, \ldots , \ell_h)$
be a non-increasing maximal chain for the prefix order ($\PMCI$)
in $\CFL_{in}(w)$, where $\ell_1$ is the first element in
$\CFL_{in}(w)$.
If $\ell_{1} \cdots  \ell_h$ is not an inverse Lyndon word, then
$h > 2$ and there are $i$, $j$, $1 \leq i < j < h$, such that
\begin{itemize}
\item[(1)]
$\ell_{1} = \ell_{2} = \ldots =  \ell_i \not = \ell_{i + 1}$,
\item[(2)]
$\ell_{i + 1} \cdots  \ell_j$ is a prefix of $\ell_1$
\item[(3)]
$\ell_{i + 1} \cdots  \ell_{j +1}$ is not
a prefix of $\ell_{1}$
\item[(4)]
We have
\begin{equation} \label{EQ5}
\ell_{1} \ll \ell_{i + 1} \cdots  \ell_{j +1}
\end{equation}
More specifically, there are words $r, s, s' \in \Sigma^*$ and $a, b \in \Sigma$, $a < b$ such that
\begin{equation} \label{EQ6}
\ell_{1} = \ell_{i + 1} \cdots  \ell_j ras, \quad
\ell_{j + 1} = rbs'
\end{equation}
\item[(5)]
We have
$$p \bar{p} =\ell_{1} \cdots  \ell_{i} \ell_{i + 1} \cdots  \ell_{j} rb$$
where $(p,\overline{p})$
is the canonical pair associated with $w$.
\end{itemize}
\end{proposition}
\begdim
Let $w \in \Sigma^+$, let $(\ell_1, \ldots , \ell_h)$
be a non-increasing maximal chain for the prefix order ($\PMCI$)
in $\CFL_{in}(w)$, where $\ell_1$ is the first element in
$\CFL_{in}(w)$.

If $\ell_{1} \cdots  \ell_h$ is not an inverse Lyndon word, then
by Proposition \ref{Algo1}, items (1)-(4) are satisfied.
Set $p_1 = \ell_{1} \cdots  \ell_{i}$, $p_2 = \ell_{i + 1} \cdots  \ell_{j} rb$.
We claim that
$$p_1p _2 =\ell_{1} \cdots  \ell_{i} \ell_{i + 1} \cdots  \ell_{j} rb = p \bar{p}$$
where $(p,\overline{p})$
is the canonical pair associated with $w$.

Indeed, by Eq. \eqref{EQ6}, we have
$$p_1 = \ell_{1} \cdots  \ell_{i} \ll \ell_{i + 1} \cdots  \ell_{j} rb = p_2$$
which implies
$$p_1p _2 \ll \ell_{i + 1} \cdots  \ell_{j} rb$$
where $\ell_{i + 1} \cdots  \ell_{j} rb$ is a suffix of $p_1p _2$, that is,
$p_1p _2$ is not an inverse Lyndon word.

Moreover, for each proper nonempty prefix $u$ of $p_2$, $u$ is also a
proper prefix of $\ell_1$, hence there is $v \in \Sigma^*$ such that
$p_1 u = (\ell_{1})^i u = (uv)^i u$ is a sesquipower of the
anti-Lyndon word $\ell_1$, hence, by Proposition \ref{InvPreN},
$p_1u$ is an inverse Lyndon word.
Moreover, by Lemma \ref{lem:inverse-Lyndon-word-prefix} and Proposition \ref{InW1},
the word $p_1$ and all its nonempty prefixes are inverse Lyndon words.
This shows that $p_1p _2$ is the shortest
prefix of $\ell_{1} \cdots  \ell_h$ which is not an inverse Lyndon word,
hence, by Proposition \ref{shortest},
$p_1p _2 = p \bar{p}$, where $(p,\overline{p})$
is the canonical pair associated with $w$.
This finishes the proof.
\enddim

Proposition \ref{step1} shows that to obtain $\ICFL(w)$ from $\CFL_{in}(w)$ we can limit ourselves to the case
in which $\CFL_{in}(w) = (\ell_1, \ldots , \ell_h)$ is a chain with respect to the prefix order.
Now, if $\ell_{1} \cdots  \ell_h$ is an inverse Lyndon word, then
$\ICFL(w)$ can be easily obtained from
$\CFL_{in}(w)$ by Proposition \ref{Algo2}.
Otherwise, Proposition \ref{Algo3} allows us to determine
the canonical pair $(p,\overline{p})$ associated with $w$
from $\CFL_{in}(w)$ and then, recursively, $\ICFL(w)$.
Examples \ref{Finale1} and \ref{Finale2} should clarify this procedure.

\begin{example} \label{Finale1}
{\rm Let $\Sigma = \{a,b \}$ with $a < b$.
Let us consider $w = babaababaababab$.
The word $w$ is not an inverse Lyndon word because
$w = babaababaababab \ll babab$.
Moreover
$$\CFL_{in}(w) = (babaa, babaa, ba, ba, b) = (\ell_1, \ell_2, \ell_3, \ell_4, \ell_5)$$
Of course, $(\ell_1, \ell_2, \ell_3, \ell_4, \ell_5)$
is a non-increasing maximal chain for the prefix order
in $\CFL_{in}(w)$, where $\ell_1$ is the first element in
$\CFL_{in}(w)$.
Items (1)-(3) in Proposition \ref{Algo3} are satisfied with
$i = 2$ and $j = 4$.
As for Eq. \eqref{EQ6} in Proposition \ref{Algo3}, we have
$$\ell_1 = babaa = \ell_3 \ell_4 a, \quad \ell_5 = b$$
By item (5) in Proposition \ref{Algo3}, for the canonical pair
$(p,\overline{p})$ associated with $w$, we have
$$p \bar{p} = \ell_1 \ell_2 \ell_3 \ell_4 \ell_5$$
We consider Proposition \ref{characterization} to determine $p$ and $\bar{p}$.
By item (2) in Proposition \ref{characterization}, $\bar{p}$ is different from $b$
because $w$ does not start with $a$, hence $\bar{p}$ ends with $ab$ and we have to look for the occurrences of the
factor $aa$ in $p \overline{p}$. There are two occurrences of $aa$ as a factor of $p \bar{p}$
and, by applying Proposition \ref{characterization}, we see that
$$p = babaa babaa = \ell_1 \ell_2, \quad \bar{p} = babab = \ell_3 \ell_4 \ell_5$$
Finally, by Proposition \ref{InW1}, $\bar{p} = babab$
is an inverse Lyndon word and, by Definition \ref{def:ICFL}, we have
$$\ICFL(w) = (p,\overline{p}) = (babaa babaa, babab)$$}
\end{example}

\begin{example} \label{Finale2}
{\rm Let $\Sigma = \{a,b,c,d \}$, $a < b < c < d$,
and $y = dabadabdabdabdadac$
(see Example \ref{groupings}).
We have
$$\CFL_{in}(y) = (daba, dab, dab, dab, dadac)$$
By Proposition \ref{step1},
$$\ICFL(y) = (\ICFL(daba dab dab dab), \ICFL(dadac))$$
The word $z = dabadabdabdab$ is not an inverse Lyndon word
and, by Proposition \ref{Algo3},
for the canonical pair
$(p,\overline{p})$ associated with $w$, we have
$p \bar{p} = daba dabd$. Thus $p = daba$,
$\bar{p} = dabd$. Then, we compute
$\ICFL(dab dab dab)$. By Propositions \ref{InW1} and \ref{Algo2},
$\ICFL(dab dab dab) = ((dab)^3)$.
Therefore, $m_1 = p$ and
$$\ICFL(daba dab dab dab) = (daba, \ICFL(dab dab dab)) = (daba, (dab)^3)$$
Of course $\ICFL(dadac) = dadac$, hence
$\ICFL(y) = (daba, (dab)^3, dadac)$.}
\end{example}

\section{Conclusions} \label{Conc}

The Lyndon factorization and the canonical inverse Lyndon factorization play a crucial role in various applications of sequence comparison and combinatorial pattern matching
\cite{BKKP19,DLT22,IS22}.
Although the two types of factorization are related by the notion of the anti-Lyndon word, their connections still remained unexplored.
This paper addresses this open problem.
More precisely, the main contribution is to show how to obtain the Lyndon factorization $\CFL_{in}(w)$
of $w$ with respect to the inverse order from the canonical inverse Lyndon factorization $\ICFL(w)$
and vice versa how to group factors of $\CFL_{in}(w)$ to obtain $\ICFL(w)$.
This result on the connection between the classical Lyndon factorization
and the unexplored inverse Lyndon factorization opens up new research directions.
For instance, the search for new bijective variants of the Burrows Wheeler Transform, based on multisets of inverse words instead of multisets
of Lyndon words.
On the other hand the characterizations provided in the paper can be used to better explore whether many results already known for $\CFL(w)$ can be easily extended to $\ICFL(w)$.
Finally, an interesting question is whether it is possible to find an alternative definition of the canonical inverse Lyndon factorization based directly only on the grouping of $\CFL_{in}(w)$. The current definition is based on the use of the pair $(p, \bar{p})$ and is not immediate to understand.

\section*{Acknowledgments}
This research was supported by the grant MIUR 2022YRB97K, PINC, Pangenome Informatics: from Theory to Applications
and by INdAM-GNCS Project 2023.

\bibliographystyle{unsrt}
\bibliography{inverse2023.bib}

\begin{thebibliography}{10}

\bibitem{inverseLyndon}
Paola Bonizzoni, C.~De~Felice, Rocco Zaccagnino, and Rosalba Zizza.
\newblock Inverse {L}yndon words and inverse {L}yndon factorizations of words.
\newblock {\em Adv. Appl. Math.}, 101:281--319, 2018.

\bibitem{Lyndon}
Kuo-Tsai Chen, Ralph~H. Fox, and Roger~C. Lyndon.
\newblock Free {D}ifferential calculus, {IV. T}he {Q}uotient {G}roups of the {L}ower {C}entral {S}eries.
\newblock {\em Ann. Math.}, 68:81--95, 1958.

\bibitem{duval}
Jean{-}Pierre Duval.
\newblock Factorizing {W}ords over an {O}rdered {A}lphabet.
\newblock {\em J. Algorithms}, 4(4):363--381, 1983.

\bibitem{apostolico-crochemore-1995}
Alberto Apostolico and Maxime Crochemore.
\newblock Fast parallel {L}yndon factorization with applications.
\newblock {\em Mathematical Systems Theory}, 28(2):89--108, 1995.

\bibitem{PerRes}
Dominique Perrin and Antonio Restivo.
\newblock Enumerative combinatorics on words.
\newblock In Miklos Bona, editor, {\em Handbook of Enumerative Combinatorics}. {CRC} Press, 2015.

\bibitem{BKKP19}
Hideo Bannai, Juha K{\"{a}}rkk{\"{a}}inen, Dominik K{\"{o}}ppl, and Marcin Piatkowski.
\newblock Indexing the bijective {BWT}.
\newblock In Nadia Pisanti and Solon~P. Pissis, editors, {\em 30th Annual Symposium on Combinatorial Pattern Matching, {CPM} 2019, June 18-20, 2019, Pisa, Italy}, volume 128 of {\em LIPIcs}, pages 17:1--17:14. Schloss Dagstuhl - Leibniz-Zentrum f{\"{u}}r Informatik, 2019.

\bibitem{GilScot}
Joseph~Yossi Gil and David~Allen Scott.
\newblock A {B}ijective {S}tring {S}orting {T}ransform.
\newblock {\em CoRR}, abs/1201.3077, 2012.

\bibitem{varBWT}
Manfred Kufleitner.
\newblock On {B}ijective {V}ariants of the {B}urrows-{W}heeler {T}ransform.
\newblock In {\em Proceedings of the Prague Stringology Conference 2009, Prague, Czech Republic, August 31 - September 2, 2009}, pages 65--79, 2009.

\bibitem{GeReu}
Ira~M. Gessel and Christophe Reutenauer.
\newblock Counting {P}ermutations with {G}iven {C}ycle {S}tructure and {D}escent {S}et.
\newblock {\em J. Comb. Theory, Ser. {A}}, 64(2):189--215, 1993.

\bibitem{GeResReu}
Ira~M. Gessel, Antonio Restivo, and Christophe Reutenauer.
\newblock A bijection between words and multisets of necklaces.
\newblock {\em Eur. J. Comb.}, 33(7):1537--1546, 2012.

\bibitem{restivo-sorting-2014}
Sabrina Mantaci, Antonio Restivo, Giovanna Rosone, and Marinella Sciortino.
\newblock Suffix array and {L}yndon factorization of a text.
\newblock {\em J. Discrete Algorithms}, 28:2--8, 2014.

\bibitem{lyndonAcc}
Nico Bertram, Jonas Ellert, and Johannes Fischer.
\newblock Lyndon words accelerate suffix sorting.
\newblock In Petra Mutzel, Rasmus Pagh, and Grzegorz Herman, editors, {\em 29th Annual European Symposium on Algorithms, {ESA} 2021, September 6-8, 2021, Lisbon, Portugal (Virtual Conference)}, volume 204 of {\em LIPIcs}, pages 15:1--15:13. Schloss Dagstuhl - Leibniz-Zentrum f{\"{u}}r Informatik, 2021.

\bibitem{Nyldon2}
{\'{E}}milie Charlier, Manon Philibert, and Manon Stipulanti.
\newblock Nyldon words.
\newblock {\em J. Comb. Theory, Ser. {A}}, 167:60--90, 2019.

\bibitem{DRR}
Francesco Dolce, Antonio Restivo, and Christophe Reutenauer.
\newblock On generalized {L}yndon words.
\newblock {\em Theor. Comput. Sci.}, 777:232--242, 2019.

\bibitem{PZ20}
Mickael Postic and Luca~Q. Zamboni.
\newblock \emph{{\(\omega\)}}-{L}yndon words.
\newblock {\em Theor. Comput. Sci.}, 809:39--44, 2020.

\bibitem{BW20}
Amanda Burcroff and Eric Winsor.
\newblock Generalized {L}yndon factorizations of infinite words.
\newblock {\em Theor. Comput. Sci.}, 809:30--38, 2020.

\bibitem{NyldonNew}
Swapnil Garg.
\newblock New results on nyldon words and nyldon-like sets.
\newblock {\em Advances in Applied Mathematics}, 131:102249, 2021.

\bibitem{inverseLyndonTCS21}
Paola Bonizzoni, C.~De~Felice, Rocco Zaccagnino, and Rosalba Zizza.
\newblock On the longest common prefix of suffixes in an inverse {L}yndon factorization and other properties.
\newblock {\em Theor. Comput. Sci.}, 862:24--41, 2021.

\bibitem{bpr}
Jean Berstel, Dominique Perrin, and Christophe Reutenauer.
\newblock {\em {C}odes and {A}utomata}.
\newblock Encyclopedia of Mathematics and its Applications 129, Cambridge University Press, 2009.

\bibitem{CK}
Christian Choffrut and Juhani Karhum\"{a}ki.
\newblock Combinatorics of {W}ords.
\newblock In Grzegorz Rozenberg and Arto Salomaa, editors, {\em Handbook of Formal Languages, Vol. 1}, pages 329--438. Springer-Verlag, Berlin, Heidelberg, 1997.

\bibitem{Lo}
M.~Lothaire.
\newblock {\em Algebraic {C}ombinatorics on {W}ords, {E}ncyclopedia {M}ath. {A}ppl.}, volume~90.
\newblock Cambridge University Press, 1997.

\bibitem{lothaire}
M.~Lothaire.
\newblock {\em Applied {C}ombinatorics on {W}ords}.
\newblock Cambridge University Press, 2005.

\bibitem{reu}
Christophe Reutenauer.
\newblock Free {L}ie algebras.
\newblock In {\em Handbook of {A}lgebra, London Mathematical Society Monographs}. Oxford Science Publications, 1993.

\bibitem{CHL07}
Maxime Crochemore, Christophe Hancart, and Thierry Lecroq.
\newblock {\em Algorithms on strings}.
\newblock Cambridge University Press, 2007.

\bibitem{Bannai15}
Hideo Bannai, I~Tomohiro, Shunsuke Inenaga, Yuto Nakashima, Masayuki Takeda, and Kazuya Tsuruta.
\newblock A new characterization of maximal repetitions by {L}yndon trees.
\newblock In {\em Proceedings of the Twenty-Sixth Annual {ACM-SIAM} Symposium on Discrete Algorithms, {SODA} 2015, San Diego, CA, USA, January 4-6, 2015}, pages 562--571, 2015.

\bibitem{Knuth}
Donald~E. Knuth.
\newblock {\em The {A}rt of {C}omputer {P}rogramming, {V}olume 4a: {C}ombinatorial {A}lgorithms}.
\newblock Addison Wesley Longman Publishing Co., Inc., 2012.

\bibitem{FM}
Harold Fredricksen and James Maiorana.
\newblock Necklaces of beads in $k$ colors and $k$-ary de {B}rujin sequences.
\newblock {\em Discrete Math.}, 23(3):207--210, 1978.

\bibitem{GGT}
Sukhpal~Singh Ghuman, Emanuele Giaquinta, and Jorma Tarhio.
\newblock Alternative {A}lgorithms for {L}yndon {F}actorization.
\newblock In {\em Proceedings of the Prague Stringology Conference 2014, Prague, Czech Republic, September 1-3, 2014}, pages 169--178, 2014.

\bibitem{lata2020}
Paola Bonizzoni, C.~De~Felice, Rocco Zaccagnino, and Rosalba Zizza.
\newblock Lyndon words versus inverse {L}yndon words: Queries on suffixes and bordered words.
\newblock In Alberto Leporati, Carlos Mart{\'{\i}}n{-}Vide, Dana Shapira, and Claudio Zandron, editors, {\em Language and Automata Theory and Applications - 14th International Conference, {LATA} 2020, Milan, Italy, March 4-6, 2020, Proceedings}, volume 12038 of {\em Lecture Notes in Computer Science}, pages 385--396. Springer, 2020.

\bibitem{antiLyndon}
Daniele~A. Gewurz and Francesca Merola.
\newblock Numeration and enumeration.
\newblock {\em Eur. J. Comb.}, 33(7):1547--1556, 2012.

\bibitem{DLT22}
Paola Bonizzoni, C.~De~Felice, Yuri Pirola, Raffaella Rizzi, Rocco Zaccagnino, and Rosalba Zizza.
\newblock Can formal languages help pangenomics to represent and analyze multiple genomes?
\newblock In Volker Diekert and Mikhail~V. Volkov, editors, {\em Developments in Language Theory - 26th International Conference, {DLT} 2022, Tampa, FL, USA, May 9-13, 2022, Proceedings}, volume 13257 of {\em Lecture Notes in Computer Science}, pages 3--12. Springer, 2022.

\bibitem{IS22}
Paola Bonizzoni, Matteo Costantini, C.~De~Felice, Alessia Petescia, Yuri Pirola, Marco Previtali, Raffaella Rizzi, Jens Stoye, Rocco Zaccagnino, and Rosalba Zizza.
\newblock Numeric lyndon-based feature embedding of sequencing reads for machine learning approaches.
\newblock {\em Inf. Sci.}, 607:458--476, 2022.

\end{thebibliography}

\end{document}